\documentclass[11pt,a4paper]{article}
\usepackage[margin=25mm]{geometry}
\usepackage[T1]{fontenc}
\usepackage{lmodern,amsmath,amssymb,amsthm,mathtools,bm}
\usepackage{graphicx,booktabs,array,longtable,microtype,xcolor}
\usepackage[hidelinks]{hyperref}
\usepackage[nameinlink,noabbrev]{cleveref}
\usepackage{enumitem}
\usepackage{caption}
\usepackage{placeins}
\usepackage{float}

\newtheorem{theorem}{Theorem}[section]
\newtheorem{lemma}[theorem]{Lemma}
\newtheorem{proposition}[theorem]{Proposition}
\newtheorem{corollary}[theorem]{Corollary}
\theoremstyle{remark}
\newcommand{\R}{\mathbb R}
\newcommand{\C}{\mathbb C}

\newcommand{\e}{\mathbf e}
\newcommand{\Ca}{C_{\mathrm a}}
\newcommand{\Cq}{C_{\mathrm q}}
\newcommand{\Cr}{C_{\mathrm r}}

\newcommand{\code}[1]{\texttt{#1}}
\title{An L-Stable Sequential Implicit Two-Stage\\
Fourth-Order Time Discretization\\
with Second-Order Stiff Damping}
\author{Zhixin Huo\\[3pt]
\small School of Mathematics and Information Science, Henan Polytechnic University,\\
\small Jiaozuo, Henan 454003, China\\
\small School of Mechatronical Engineering, Beijing Institute of Technology,\\
\small Beijing 100081, China\\
\small\href{mailto:zhixinhuo@hpu.edu.cn}{zhixinhuo@hpu.edu.cn}}
\date{}
\begin{document}
\maketitle
\begin{abstract}

We construct a fully implicit, sequential two-stage time discretization for stiff systems using physical quantities and their rates of change. The undetermined coefficients are determined by direct Taylor matching, and the convergence analysis proves that the scheme achieves fourth-order accuracy at both the midpoint and the endpoint for fixed smooth systems. The analysis of the poles and the imaginary-axis boundary, together with the maximum modulus principle, characterizes the \(L\)-stable family. Within this family, \(C_q=5/183\) uniquely cancels the leading inverse-power term and endows the full-step amplification factor with quadratic stiff decay, while \(C_a=\frac{25-\sqrt{105}}{780}\) is the parameter choice yielding the highest accuracy. The new method simultaneously possesses the advantages of \(L\)-stability with second-order stiff damping, two-stage fourth-order accuracy, and sequential implicitness. In contrast, the fourth-order Gauss--Legendre scheme and the symmetric Hermite--Obreshkov scheme are both \(A\)-stable rather than \(L\)-stable; their variants have to sacrifice the accuracy attainable at each stage in order to improve stiff stability. For example, among \(L\)-stable methods, two-stage Radau IIA and TR--BDF2 have third- and second-order accuracy, respectively, while the SDIRK method requires five stages to achieve fourth-order accuracy. In addition, the stages of the Gauss and Radau formulations must be solved in a coupled manner. Stiff ODE tests and scalar tests verify fourth-order convergence and quadratic stiff decay, respectively. Finite-volume tests show that the method exhibits fourth-order joint convergence on smooth one-dimensional problems and fast-mode damping in two-dimensional cases. At the same acoustic error threshold, the method uses eight implicit stages, whereas three-stage Radau IIA uses twelve. Finally, the sequential formulation keeps each Newton system at the state dimension, thereby relatively reducing its size, improving computational efficiency, and facilitating parallel computation.

\end{abstract}
\noindent\textbf{Keywords:} two-stage fourth-order methods; implicit time discretization; stiff damping; L-stability; sequential advancement.
\clearpage

\section{Introduction}
To achieve consistent high-order accuracy in space and time, the accuracy of time advancement needs to match that of spatial discretization \cite{jiangshu,cockburnshu}. Common high-order time discretizations mainly include Runge--Kutta and Lax--Wendroff-type methods. The former achieve high order by advancing physical quantities through several stages within one time step. Although simple to implement, they generally repeat spatial reconstruction and flux evaluation in high-order flow solvers, increasing computational work and extending the effective stencil of the completed step.The latter use the classical Cauchy--Kovalevskaya procedure, recursively applying the governing equations to convert higher temporal derivatives into spatial derivatives. For sufficiently smooth solutions, a single-stage construction can achieve arbitrarily high formal order \cite{lax,lireview}. However, repeated differentiation is difficult to interpret classically in flows with poor regularity, and the increasing complexity of flux calculations limits its practical use in complex flows \cite{lireview}.

To balance these advantages and disadvantages, my doctoral advisor, Professor Jiequan Li, and his collaborators proposed and developed the two-stage fourth-order space--time coupled method for Lax--Wendroff-type solvers, including GRP and GKS \cite{lireview,lispacetime,lidu,gks4,pan3dgks}. The method is constructed jointly from physical quantities and their rates of change, a pair of interconnected quantities with clear physical meanings. This construction recalls the joint use of position and momentum in symplectic methods: in this analogy, the physical quantity corresponds to position and serves as a solution element, while its rate corresponds to momentum and serves as a dynamics element \cite{lireview,fengqin}. Just as symplectic methods preserve the symplectic structure of Hamiltonian systems, the two-stage fourth-order method can preserve the flow's space--time correlation through compatible local evolution and time-integrated fluxes \cite{lireview,lispacetime,fluxregularity}. Relative to Runge--Kutta constructions based on solution quantities alone, compatible flow implementations offer substantial gains in computational efficiency \cite{lidu,chenggrp,gks4} and compactness \cite{hweno4,jicompact}; suitable variable-weight extensions also improve numerical stability \cite{yuantang}. Relative to classical Lax--Wendroff constructions using higher temporal derivatives, the method greatly simplifies flux evaluation while retaining a physical interpretation of local evolution in complex flows \cite{lireview}.

For $u'=\mathcal{L}(u)$ and $\mathcal{H}=\mathcal{L}_u\mathcal{L}$, the original two-stage fourth-order time discretization reads
\begin{align*}\label{eq:liduoriginal}
 u^{n+1/2}&=u^n+\frac h2\mathcal{L}(u^n)+\frac{h^2}{8}\mathcal{H}(u^n),\notag\\
 u^{n+1}&=u^n+h\mathcal{L}(u^n)+\frac{h^2}{6}\{\mathcal{H}(u^n)+2\mathcal{H}(u^{n+1/2})\}.
\end{align*}
This original space--time coupled construction uses explicit advancement and is suited to convection-dominated nonstiff or weakly stiff problems. Practical multiphysics systems often contain rapid reactions, strong relaxation, and diffusion, and exhibit strongly stiff multiscale dynamics. The CFL stability conditions of this explicit scheme can then constrain the time step to the smallest characteristic time scale, substantially reducing efficiency. Thus professor Li explicitly emphasized the urgent need to extend the two-stage fourth-order time-space coupled method from explicit to implicit advancement \cite{lispacetime}.

Implicit integration using additional derivative information has a substantial history. Hairer and Wanner \cite{hw1973} studied general multistep--multistage--multiderivative representations and order analysis; Enright \cite{enright} and Cash \cite{cash} developed second-derivative multistep formulas and second-derivative extended backward differentiation formulas for stiff systems. Hermite-type one-step methods were studied by N\o rsett \cite{norsett}. Mazzia and Sestini investigated Hermite--Obreshkov methods and continuous spline extensions \cite{mazzia2018}, and revised the related conjugate-symplecticity conclusion in a subsequent correction \cite{mazzia2019}; Iavernaro and Mazzia \cite{iavernaro} further studied related fourth-order symplectic and conjugate-symplectic constructions. For PDEs, Tsai et al.\ \cite{tsaicw2014} studied spatial discretization of two-derivative information and applied a two-stage fourth-order semi-implicit method to the heat equation. Zeifang and Sch\"utz \cite{zeifang} combined implicit two-derivative deferred correction with DG discretization, with subsequent work addressing time parallelism and Newton adaptivity \cite{zeifangamc}. Chouchoulis and Sch\"utz \cite{chouchoulis} investigated Jacobian-free implicit multiderivative constructions and solvers. These studies show that additional evolution information supports high-order implicit integration, while practical efficiency also depends on derivative construction and implicit solution costs.

Following Li and Du's \cite{lidu} construction from physical quantities and their rates of change through undetermined coefficients and Taylor expansion, this paper independently develops a sequential implicit two-stage fourth-order time discretization. Two successively solved implicit stage equations are introduced, and Taylor coefficient matching determines the conditions for fourth-order accuracy at both the intermediate node and the full-step endpoint. Global fourth-order convergence is proved under appropriate smoothness assumptions. Characteristic analysis, pole analysis, and imaginary-axis boundary analysis, together with the maximum modulus principle, yield the L-stable relation $D=-C$ and its admissible interval. The subsequent Huo–Su study \cite{huosu} identifies, by cancelling the leading inverse-power term in the complete amplification factor, the unique choice $C=5/183$ giving quadratic stiff damping within this family. At equal fourth-order accuracy, the present scheme possesses stronger stiff-damping capability, whereas the two-stage Gauss–Legendre method and the classical symmetric fourth-order Hermite–Obreshkov method are A-stable but not L-stable. Among the schemes that are all L-stable, the present scheme attains fourth-order accuracy with two-stage advancement, whereas two-stage Radau IIA and TR–BDF2 achieve only third- and second-order accuracy, respectively \cite{radau,hairerbook,hosea}, and the SDIRK method requires five stages \cite{alexander,kennedy} to reach fourth-order accuracy. In addition, the two stages of the present scheme are solved sequentially, avoiding the coupled computation among the nonlinear stage equations of the Gauss/Radau methods.

Although the first stage of the present scheme is algebraically identical to the classical symmetric fourth-order Hermite–Obreshkov rule over a half time step, for the characterized L-stable family the complete amplification factor of the present scheme tends to zero in the deep-stiff limit, which distinguishes it from the classical symmetric Hermite full step and its simple composition of two half steps.
Hairer–Wanner's general theory \cite{hw1973} provides an analytical language for several classes of integration methods. Expressing a method in a unified notation does not by itself establish that its particular sequential structure, coefficient set, and stability conclusions have already been supplied by that general representation. Here the specific implicit equations are developed through undetermined coefficients and Taylor matching; standard multiderivative representations and order conditions are used only for subsequent independent verification and comparison, and are not the starting point of the coefficient derivation or the fourth-order proof.

The paper details Newton iteration and the one- and two-dimensional implementations after spatial discretization: a Lax–Wendroff solver is used to obtain the physical quantities and their rates of change at the interfaces, and WENO spatial discretization is combined with the implicit two-stage fourth-order time discretization to form a space–time coupled numerical scheme. Numerical verification begins with stiff ODEs, including a linear stiff system, Robertson kinetics, the eight-component HIRES reaction problem, and the van der Pol oscillator, to examine classical fourth-order convergence, L-stable attenuation, and the predicted quadratic stiff decay. Joint space–time refinement yields fourth-order convergence for smooth one-dimensional advection and Burgers problems, and unsplit two-dimensional transport–relaxation tests demonstrate fast-mode damping. At matched error targets, the method requires fewer implicit stages than selected L-stable baseline methods. The two stages are solved sequentially, which reduces the size of the individual Newton system compared with the coupled Radau IIA formulation.

\numberwithin{equation}{section}
\section{Construction and fourth-order accuracy by Taylor matching}\label{sec:method}
First, the exact midpoint and endpoint values are obtained from the integral evolution, the trajectory chain rule, and the integrated expansion of the evolution operator. Next, two numerical equations with undetermined coefficients, constructed from the physical quantities and their rates of change, are introduced and expanded directly in powers of the time step via Taylor expansion. Comparing the numerical coefficients with those of the corresponding exact values determines the method. The calculation below retains the distinct vector contractions. No Runge–Kutta order conditions are used to determine the coefficients.

\subsection{Integral evolution and exact midpoint and endpoint expansions}
Consider the autonomous evolution problem
\begin{equation}\label{eq:ode}
 \frac{du}{dt}=\mathcal{L}(u),\qquad
 u(0)=u_0\in\R^N,
\end{equation}
where $u=(u_1,\ldots,u_N)^T$ is the vector of physical quantities and the operator $\mathcal{L}:\R^N\to\R^N$ is sufficiently smooth. Throughout this paper, the evolution operator depends on the state only. Spatial coefficients, source laws, and prescribed boundary data are independent of time; the solution and its trajectory derivatives evolve with time. With $h=\Delta t>0$, its exact integral form is
\begin{equation}\label{eq:integraloriginal}
 u(t_n+\theta h)=u(t_n)+\int_{t_n}^{t_n+\theta h}\mathcal{L}(u(t))\,dt,
 \qquad 0\le\theta\le1.
\end{equation}
To avoid repeated arguments, an ordinary pair of parentheses groups an operator expression with a common evaluation state; for example,
$(\mathcal L_u\mathcal L)(u^n)=\mathcal L_u(u^n)\mathcal L(u^n)$.
The trailing state argument applies to the operators and their state derivatives, while displayed directions and numerical-state differences retain their stated values. In the chain-rule identity below, all operators are evaluated along $u(t)$ and their arguments are suppressed.
Write $Q(t)=\mathcal{L}(u(t))$. Differentiation along the solution gives
\begin{equation}\label{eq:vector}
\begin{aligned}
 Q'&=\mathcal L_u\mathcal L,\\
 Q''&=\mathcal L_{uu}[\mathcal L,\mathcal L]+\mathcal L_u^2\mathcal L,\\
 Q'''&=\mathcal L_{uuu}[\mathcal L,\mathcal L,\mathcal L]
       +3\mathcal L_{uu}[\mathcal L,\mathcal L_u\mathcal L]+\mathcal L_u\bigl(\mathcal L_{uu}[\mathcal L,\mathcal L]\bigr)
       +\mathcal L_u^3\mathcal L.
\end{aligned}
\end{equation}
Here $\mathcal L_u$ is a Jacobian matrix, and higher state derivatives are multilinear maps; powers of $\mathcal L_u$ mean matrix composition. In particular, the two mixed contractions in $Q'''$ are distinct for a general vector field. Their scalar counterparts coincide and have coefficients $3+1=4$. We retain them separately throughout the vector calculation.

Taylor expansion inside the integral yields
\begin{equation}\label{eq:integraltaylor}
\begin{aligned}
 &\int_{t_n}^{t_n+\theta h}\mathcal{L}(u(t))\,dt\\
 &=\int_{t_n}^{t_n+\theta h}
 \left[Q_n+(t-t_n)Q'_n+\frac{(t-t_n)^2}{2!}Q''_n
             +\frac{(t-t_n)^3}{3!}Q'''_n+O((t-t_n)^4)\right]dt\\
 &=\theta hQ_n+\frac{(\theta h)^2}{2!}Q'_n
             +\frac{(\theta h)^3}{3!}Q''_n
             +\frac{(\theta h)^4}{4!}Q'''_n+O(h^5).
\end{aligned}
\end{equation}
In the following two formulas, the common exact state $u(t_n)$ is written once after the closing parenthesis. The integral relation \eqref{eq:integraloriginal}, the trajectory derivatives \eqref{eq:vector}, and the evaluated integral \eqref{eq:integraltaylor} jointly determine the exact nodal values. Substituting $\theta=1/2$ into those three relations gives
\begin{equation}\label{eq:exactmid}
\begin{aligned}
&u(t_n+h/2)={}u(t_n)+\Bigl(\frac h2\mathcal L+\frac{h^2}{8}\mathcal L_u\mathcal L+\frac{h^3}{48}\bigl(\mathcal L_{uu}[\mathcal L,\mathcal L]+\mathcal L_u^2\mathcal L\bigr)\\
&\quad\quad\quad+\frac{h^4}{384}\bigl(\mathcal L_{uuu}[\mathcal L,\mathcal L,\mathcal L]
       +3\mathcal L_{uu}[\mathcal L,\mathcal L_u\mathcal L]+\mathcal L_u\bigl(\mathcal L_{uu}[\mathcal L,\mathcal L]\bigr)+\mathcal L_u^3\mathcal L\bigr)\Bigr)(u(t_n))+O(h^5).
\end{aligned}
\end{equation}
Substituting $\theta=1$ into the same three relations gives
\begin{equation}\label{eq:exactend}
\begin{aligned}
&u(t_n+h)={}u(t_n)+\Bigl(h\mathcal L+\frac{h^2}{2}\mathcal L_u\mathcal L+\frac{h^3}{6}\bigl(\mathcal L_{uu}[\mathcal L,\mathcal L]+\mathcal L_u^2\mathcal L\bigr)\\
&\quad\quad\quad+\frac{h^4}{24}\bigl(\mathcal L_{uuu}[\mathcal L,\mathcal L,\mathcal L]
       +3\mathcal L_{uu}[\mathcal L,\mathcal L_u\mathcal L]+\mathcal L_u\bigl(\mathcal L_{uu}[\mathcal L,\mathcal L]\bigr)+\mathcal L_u^3\mathcal L\bigr)\Bigr)(u(t_n))+O(h^5).
\end{aligned}
\end{equation}

The two formulas above are the half-step and full-step exact values obtained by integrating the evolution law; the Taylor expansion used in this calculation is the operator expansion inside the integral in \eqref{eq:integraltaylor}. They provide the two distinct references for the coefficient comparison below. Exact values are written as $u(t_n+\theta h)$, whereas $u^n$, $u^{n+1/2}$, and $u^{n+1}$ denote numerical states.

\subsection{A sequential implicit two-stage ansatze}
For a single step, $u^n$ is known, while $u^{n+1/2}$ and $u^{n+1}$ approximate the exact values at $t_n+h/2$ and $t_n+h$. Denote the trajectory derivative of the evolution operator by $\mathcal{H}(u)=\frac{d}{dt}\mathcal{L}(u(t))\big|_{u(t)=u}=\mathcal{L}_u(u)\mathcal{L}(u)$. The first implicit stage is
\begin{equation}\label{eq:firstansatz}
\begin{aligned}
 u^{n+1/2}={}&u^n+h\left[c_0\mathcal{L}(u^n)+c_{1/2}\mathcal{L}(u^{n+1/2})\right]+h^2\left[d_0\mathcal{H}(u^n)+d_{1/2}\mathcal{H}(u^{n+1/2})\right].
\end{aligned}
\end{equation}
The second implicit stage uses the initial state, the computed intermediate state, and the endpoint:
\begin{equation}\label{eq:secondansatz}
\begin{aligned}
 u^{n+1}={}&u^n+h\left[a_0\mathcal{L}(u^n)+a_{1/2}\mathcal{L}(u^{n+1/2})+a_1\mathcal{L}(u^{n+1})\right]\\
 &\quad\;+h^2\left[b_0\mathcal{H}(u^n)+b_{1/2}\mathcal{H}(u^{n+1/2})+b_1\mathcal{H}(u^{n+1})\right].
\end{aligned}
\end{equation}
All ten coefficients are initially undetermined. The first equation is solved before the second; no simultaneous nonlinear solve for the two unknown stages is introduced. We impose fourth-order accuracy at each of the two nodes.

\subsection{Direct Taylor expansion of the numerical stage equations}
For the local comparison, take $u^n=u(t_n)$. The exact solution at $t_n+\theta h$, $\theta\in\{1/2,1\}$ has the Taylor expansion
\begin{equation}\label{eq:u_taylor_before_substitution}
\begin{aligned}
 u(t_n+\theta h)={}&u(t_n)+\theta h u_t(t_n)+\frac{(\theta h)^2}{2!}u_{tt}(t_n)
 +\frac{(\theta h)^3}{3!}u_{ttt}(t_n)+\frac{(\theta h)^4}{4!}u_{tttt}(t_n)+O(h^5),
\end{aligned}
\end{equation}
where the following trajectory identities are evaluated at $t_n$:
\begin{equation}\label{eq:trajectory_derivatives}
\begin{aligned}
 u_t&=\mathcal L,\qquad u_{tt}=\mathcal L_u\mathcal L,\qquad
 u_{ttt}=\mathcal L_{uu}[\mathcal L,\mathcal L]+\mathcal L_u^2\mathcal L,\\
 u_{tttt}&=\mathcal L_{uuu}[\mathcal L,\mathcal L,\mathcal L]
 +3\mathcal L_{uu}[\mathcal L,\mathcal L_u\mathcal L]
 +\mathcal L_u\bigl(\mathcal L_{uu}[\mathcal L,\mathcal L]\bigr)+\mathcal L_u^3\mathcal L.
\end{aligned}
\end{equation}
%Using $u_t=\mathcal L(u)$ and the trajectory identities \eqref{eq:vector} gives the corresponding result entirely in terms of $\mathcal L$:
Let $u^{n+\theta}=u(t_n+\theta h)$, $u^n=u(t_n)$, then we have
\begin{equation}\label{eq:u_taylor_L_result}
\begin{aligned}
&u^{n+\theta}={}u^n+\Bigl(\theta h\mathcal L+\frac{(\theta h)^2}{2}\mathcal L_u\mathcal L+\frac{(\theta h)^3}{6}\bigl(\mathcal L_{uu}[\mathcal L,\mathcal L]+\mathcal L_u^2\mathcal L\bigr)\\
 &\quad\quad\quad+\frac{(\theta h)^4}{24}\bigl(\mathcal L_{uuu}[\mathcal L,\mathcal L,\mathcal L]
       +3\mathcal L_{uu}[\mathcal L,\mathcal L_u\mathcal L]+\mathcal L_u\bigl(\mathcal L_{uu}[\mathcal L,\mathcal L]\bigr)+\mathcal L_u^3\mathcal L\bigr)\Bigr)(u^n)+O(h^5).
\end{aligned}
\end{equation}
In \eqref{eq:exactmid} and \eqref{eq:exactend}, $u(t_n+\theta h)$ denotes the exact solution; $u^{n+\theta}$ denotes a numerical stage determined by \eqref{eq:firstansatz}--\eqref{eq:secondansatz}. Expand those implicit equations and match powers of $h$ successively. At each degree, use only the numerical-state coefficients fixed at preceding degrees. Terms through degree three suffice in $h\mathcal L$, and terms through degree two suffice in $h^2\mathcal H$, this procedure determines the fourth-degree numerical coefficient. 
%without presupposing fourth-order accuracy.

For arbitrary directions $\xi,\zeta\in\R^N$, the product rule gives the following pointwise identities; their common state argument is omitted:
\begin{equation}\label{eq:H_state_derivatives}
\begin{aligned}
 \mathcal H&=\mathcal L_u\mathcal L,\\
 \mathcal H_u[\xi]&=\mathcal L_{uu}[\xi,\mathcal L]+\mathcal L_u^2\xi,\\
 \mathcal H_{uu}[\xi,\zeta]&=\mathcal L_{uuu}[\xi,\zeta,\mathcal L]
       +\mathcal L_{uu}[\xi,\mathcal L_u\zeta]+\mathcal L_{uu}[\zeta,\mathcal L_u\xi]
       +\mathcal L_u\bigl(\mathcal L_{uu}[\xi,\zeta]\bigr).
\end{aligned}
\end{equation}
Taylor-expand the operators about $u^n$. The first equalities below are state Taylor expansions. The simplified second equalities use the numerical-state coefficients already matched through degree three for $\mathcal L$ and degree two for $\mathcal H$, together with \eqref{eq:H_state_derivatives}. They apply during the successive coefficient comparison for $\theta\in\{1/2,1\}$:
\begin{equation}\label{eq:Lstate_direct}
\begin{aligned}
 &\mathcal L(u^{n+\theta})
 ={}\Bigl(\mathcal L+\mathcal L_u[u^{n+\theta}-u^n]+\frac12\mathcal L_{uu}[u^{n+\theta}-u^n,u^{n+\theta}-u^n]\\
 &\quad\quad\quad\quad+\frac16\mathcal L_{uuu}[u^{n+\theta}-u^n,u^{n+\theta}-u^n,u^{n+\theta}-u^n]\Bigr)(u^n)+O(h^4)\\
 &\quad\quad\quad\;\;={}\Bigl(\mathcal L+\theta h\mathcal L_u\mathcal L+\frac{(\theta h)^2}{2}\bigl(\mathcal L_{uu}[\mathcal L,\mathcal L]+\mathcal L_u^2\mathcal L\bigr)\\
 & \quad\quad\quad\quad+\frac{(\theta h)^3}{6}\bigl(\mathcal L_{uuu}[\mathcal L,\mathcal L,\mathcal L]
      +3\mathcal L_{uu}[\mathcal L,\mathcal L_u\mathcal L]+\mathcal L_u\bigl(\mathcal L_{uu}[\mathcal L,\mathcal L]\bigr)
       +\mathcal L_u^3\mathcal L\bigr)\Bigr)(u^n)+O(h^4),
\end{aligned}
\end{equation}
and
\begin{equation}\label{eq:Hstate_direct}
\begin{aligned}
 &\mathcal H(u^{n+\theta})
 ={}\Bigl(\mathcal H+\mathcal H_u[u^{n+\theta}-u^n]+\frac12\mathcal H_{uu}[u^{n+\theta}-u^n,u^{n+\theta}-u^n]\Bigr)(u^n)+O(h^3)\\
 &\quad\quad\quad\;\;\;={}\Bigl(\mathcal L_u\mathcal L
       +\theta h\bigl(\mathcal L_{uu}[\mathcal L,\mathcal L]+\mathcal L_u^2\mathcal L\bigr)\\
 &\quad\quad\quad\quad+\frac{(\theta h)^2}{2}\bigl(\mathcal L_{uuu}[\mathcal L,\mathcal L,\mathcal L]
       +3\mathcal L_{uu}[\mathcal L,\mathcal L_u\mathcal L]+\mathcal L_u\bigl(\mathcal L_{uu}[\mathcal L,\mathcal L]\bigr)
       +\mathcal L_u^3\mathcal L\bigr)\Bigr)(u^n)+O(h^3).
\end{aligned}
\end{equation}
The common expansion state $u^n$ is indicated once after each grouped expression. The operator value $\mathcal L$ and its Jacobian $\mathcal L_u$ remain distinct. The two mixed vector contractions retain separate factors $3$ and $1$. The successive lower-order comparisons justify each simplified operator coefficient.

Continuing with \eqref{eq:Lstate_direct}--\eqref{eq:Hstate_direct} gives the following form, read degree by degree under the preceding conditions:
\begin{equation}\label{eq:midpointansatzexpansion}
\begin{aligned}
&u^{n+1/2}={}u^n+\Bigl(h(c_0+c_{1/2})\mathcal L
       +h^2\bigl(\tfrac12c_{1/2}+d_0+d_{1/2}\bigr)\mathcal L_u\mathcal L\\
&\quad\quad+h^3\bigl(\tfrac18c_{1/2}+\tfrac12d_{1/2}\bigr)
       \bigl(\mathcal L_{uu}[\mathcal L,\mathcal L]+\mathcal L_u^2\mathcal L\bigr)+h^4\bigl(\tfrac1{48}c_{1/2}+\tfrac18d_{1/2}\bigr)\\
&\times\bigl(\mathcal L_{uuu}[\mathcal L,\mathcal L,\mathcal L]+3\mathcal L_{uu}[\mathcal L,\mathcal L_u\mathcal L]
       +\mathcal L_u\bigl(\mathcal L_{uu}[\mathcal L,\mathcal L]\bigr)
       +\mathcal L_u^3\mathcal L\bigr)\Bigr)(u^n)+O(h^5).
\end{aligned}
\end{equation}
After matching the first stage, treat the endpoint equation in the same order, retaining the computed intermediate state. Each endpoint coefficient below likewise uses the preceding endpoint matching conditions:
\begin{equation}\label{eq:endpointansatzexpansion}
\begin{aligned}
&u^{n+1}={}u^n+\Bigl(h(a_0+a_{1/2}+a_1)\mathcal L+h^2\bigl(\tfrac12a_{1/2}+a_1+b_0+b_{1/2}+b_1\bigr)\mathcal L_u\mathcal L\\
 &\quad+h^3\bigl(\tfrac18a_{1/2}+\tfrac12a_1+\tfrac12b_{1/2}+b_1\bigr)
       \bigl(\mathcal L_{uu}[\mathcal L,\mathcal L]+\mathcal L_u^2\mathcal L\bigr)+h^4\bigl(\tfrac1{48}a_{1/2}+\tfrac16a_1+\tfrac18b_{1/2}+\tfrac12b_1\bigr)\\
 &\quad\times\bigl(\mathcal L_{uuu}[\mathcal L,\mathcal L,\mathcal L]
       +3\mathcal L_{uu}[\mathcal L,\mathcal L_u\mathcal L]+\mathcal L_u\bigl(\mathcal L_{uu}[\mathcal L,\mathcal L]\bigr)
       +\mathcal L_u^3\mathcal L\bigr)\Bigr)(u^n)+O(h^5).
\end{aligned}
\end{equation}
The corresponding exact values \eqref{eq:exactmid} and \eqref{eq:exactend} now provide the two coefficient comparisons. %The numerical states have not been replaced by exact intermediate or endpoint data.

\subsection{Coefficient conditions and the completed two-stage method}
Matching the powers $h,h^2,h^3,h^4$ successively in \eqref{eq:midpointansatzexpansion} with the exact half-step \eqref{eq:exactmid} gives
\begin{equation}\label{eq:firstconditions}
\begin{aligned}
 c_0+c_{1/2}&=\frac12,\\
 \frac12c_{1/2}+d_0+d_{1/2}&=\frac18,\\
 \frac18c_{1/2}+\frac12d_{1/2}&=\frac1{48},\\
 \frac1{48}c_{1/2}+\frac18d_{1/2}&=\frac1{384}.
\end{aligned}
\end{equation}
The unique solution is
\begin{equation}\label{eq:firstcoefficients}
 c_0=c_{1/2}=\frac14,\qquad d_0=\frac1{48},\qquad d_{1/2}=-\frac1{48}.
\end{equation}
The same successive comparison of \eqref{eq:endpointansatzexpansion} with the exact endpoint \eqref{eq:exactend} gives
\begin{equation}\label{eq:secondconditions}
\begin{aligned}
 a_0+a_{1/2}+a_1&=1,\\
 \frac12a_{1/2}+a_1+b_0+b_{1/2}+b_1&=\frac12,\\
 \frac18a_{1/2}+\frac12a_1+\frac12b_{1/2}+b_1&=\frac16,\\
 \frac1{48}a_{1/2}+\frac16a_1+\frac18b_{1/2}+\frac12b_1&=\frac1{24}.
\end{aligned}
\end{equation}
With $b_0=C$ and $b_{1/2}=D$, the solution is
\begin{equation}\label{eq:coefficients}
\begin{aligned}
 a_0&=\frac16+4C+\frac D2,& b_0&=C,\\
 a_{1/2}&=\frac23-8C+2D,& b_{1/2}&=D,\\
 a_1&=\frac16+4C-\frac{5D}{2},& b_1&=\frac D2-C.
\end{aligned}
\end{equation}
The separate vector contractions in \eqref{eq:vector} do not change these coefficient solutions. 
%At the midpoint, their two fourth-order matching requirements are
%\begin{equation}\label{eq:split_mid_conditions}
% 3\left(\frac{c_{1/2}}{48}+\frac{d_{1/2}}8\right)=\frac3{384},
% \qquad
% \frac{c_{1/2}}{48}+\frac{d_{1/2}}8=\frac1{384}.
%\end{equation}
%Both are the fourth equation of \eqref{eq:firstconditions}. At the endpoint they are
%\begin{equation}\label{eq:split_end_conditions}
%\begin{aligned}
% 3\left(\frac{a_{1/2}}{48}+\frac{a_1}6+\frac{b_{1/2}}8+\frac{b_1}2\right)&=\frac3{24},\\
% \frac{a_{1/2}}{48}+\frac{a_1}6+\frac{b_{1/2}}8+\frac{b_1}2&=\frac1{24},
%\end{aligned}
%\end{equation}
%which are the same fourth equation of \eqref{eq:secondconditions}. 
%The first-stage weights and the complete $(C,D)$ family remain unchanged, thus the two distinct vector terms satisfy their own required coefficients without an additional independent constraint. Consequently, the L-stable interval and the later quadratic-decay choice derived from the complete scalar response are unchanged as well. %This conclusion follows from the displayed vector identities, rather than from scalar testing alone.
The two distinct vector terms
\(3\mathcal{L}_{uu}[\mathcal{L},\mathcal{L}_u\mathcal{L}]\)
and
\(\mathcal{L}_u\bigl(\mathcal{L}_{uu}[\mathcal{L},\mathcal{L}]\bigr)\)
each satisfy their own required coefficients and are equivalent to each other;
therefore, the first-stage weights and the complete \((C,D)\) parameter family
remain unchanged. Consequently, the L-stable interval derived from the complete
scalar response and the subsequent quadratic-decay choice also remain unchanged.

Substitution into the ansatze \eqref{eq:firstansatz}-\eqref{eq:secondansatz} yields
\begin{equation}\label{eq:stage1}
\begin{aligned}
 u^{n+1/2}={}&u^n+\frac h4\left[\mathcal{L}(u^n)+\mathcal{L}(u^{n+1/2})\right]+\frac{h^2}{48}\left[\mathcal{H}(u^n)-\mathcal{H}(u^{n+1/2})\right],
\end{aligned}
\end{equation}
and
\begin{equation}\label{eq:stage2}
\begin{aligned}
 u^{n+1}={}&u^n+h\Bigl[\left(\frac16+4C+\frac D2\right)\mathcal{L}(u^n)
                +\left(\frac23-8C+2D\right)\mathcal{L}(u^{n+1/2})+\left(\frac16+4C-\frac{5D}{2}\right)\mathcal{L}(u^{n+1})\Bigr]\\
 &+h^2\left[C \mathcal{H}(u^n)+D \mathcal{H}(u^{n+1/2})+\left(\frac D2-C\right)\mathcal{H}(u^{n+1})\right].
\end{aligned}
\end{equation}
First solve for $u^{n+1/2}$, then for $u^{n+1}$. Each unknown stage has dimension $N$. Taylor coefficient matching gives fourth-order consistency at both nodes; the local and global errors are established below.
% The exact-arithmetic scripts supplied with the code independently verify the numerical expansions and coefficient solutions.

\subsection{Taylor accuracy and global fourth-order convergence}\label{sec:order}
Assume that the exact solution exists on $[0,T]$ and $\mathcal L\in C^5$ with bounded derivatives near its compact trajectory; then $u\in C^6$ and $\mathcal H\in C^4$. Fix finite $C,D$ in \eqref{eq:coefficients}. Denote the first-stage and completed-step maps by $\Psi_h,\Phi_h$. Constants are uniform in $n,h$ for this fixed system, parameters and interval. All results concern sufficiently small $h$ and the locally consistent branches.

\begin{lemma}[Local solvability and local errors]\label{lem:defect}\label{lem:solvability}
The two stages have locally unique consistent branches. Exact starting data give the uniform bounds
\begin{equation}\label{eq:midlocal}
 \|\Psi_h(u(t_n))-u(t_n+h/2)\|\le K_{\mathrm{loc}}h^5,
\end{equation}
\begin{equation}\label{eq:local}
 \|\Phi_h(u(t_n))-u(t_n+h)\|\le K_{\mathrm{loc}}h^5.
\end{equation}
For $D=-C$, the sharper endpoint expansion is
\begin{equation}\label{eq:principal}
 \Phi_h(u(t_n))-u(t_n+h)
 =\frac{1-75C}{2880}h^5u^{(5)}(t_n)+O(h^6).
\end{equation}
\end{lemma}
\begin{proof}
The residual Jacobians are $I+O(h)$ uniformly. The implicit function theorem and inverse bound in Appendix~\ref{app:uniform_branches} give the local branches. Taylor-matched residuals $O(h^5)$ then give \eqref{eq:midlocal}--\eqref{eq:local}; the computed-midpoint substitution changes the endpoint only at $O(h^6)$. Appendix~\ref{app:local_errors} derives the uniform bounds and the sharper vector quadrature expansion \eqref{eq:principal}.
\end{proof}

\begin{theorem}[Global fourth-order convergence at both nodes]\label{thm:convergence}\label{cor:bothnodes}
Under these assumptions, let $\|u^0-u(0)\|=O(h^4)$. For sufficiently small $h$ there is a constant $K_T$ such that
\begin{equation}\label{eq:global}
 \max_{0\le nh\le T}\|u^n-u(t_n)\|\le K_T h^4,
\end{equation}
and
\begin{equation}\label{eq:bothnodesglobal}
 \max_{t_n+h\le T}\Bigl(\|u^{n+1/2}-u(t_n+h/2)\|
                   +\|u^{n+1}-u(t_n+h)\|\Bigr)\le K_T h^4.
\end{equation}
These conclusions also hold for inexact stage solves in the same neighborhood when each unscaled stage residual is at most $K_s h^5$, uniformly over the steps.
\end{theorem}
\begin{proof}
The exact sensitivity identities \eqref{eq:app_mid_sensitivity}--\eqref{eq:app_end_sensitivity}, including the midpoint chain contribution, give
\[
 \|\mathrm D\Psi_h\|\le1+K_{\mathrm{prop}}h,
 \qquad \|\mathrm D\Phi_h\|\le1+K_{\mathrm{prop}}h.
\]
With $e_n=u^n-u(t_n)$, the mean-value estimate and \cref{lem:defect} yield
\[
 \|e_{n+1}\|\le(1+K_{\mathrm{prop}}h)\|e_n\|+K_{\mathrm{loc}}h^5,
 \qquad
 \|e_n\|\le e^{K_{\mathrm{prop}}T}
             \bigl(\|e_0\|+K_{\mathrm{loc}}Th^4\bigr),
\]
which proves \eqref{eq:global}. The recomputed midpoint satisfies
\[
\begin{aligned}
 \|\Psi_h(u^n)-u(t_n+h/2)\|
 &\le\|\Psi_h(u^n)-\Psi_h(u(t_n))\|+K_{\mathrm{loc}}h^5\\
 &\le(1+K_{\mathrm{prop}}h)\|e_n\|+K_{\mathrm{loc}}h^5.
\end{aligned}
\]
This proves \eqref{eq:bothnodesglobal}. Stage residuals $O(h^5)$ add $O(h^5)$ to the recurrence by Appendix~\ref{app:inexact_newton}. The uniform step choice and induction in Appendix~\ref{app:uniform_branches} keep every input and stage inside the selected neighborhoods, closing the argument for exact and inexact solves.
\end{proof}

The theorem also applies to polynomial reaction systems with a compact invariant concentration set, including the Robertson problem; Appendix~\ref{app:robertson_regular} verifies the regularity hypotheses for that example.
Constants may depend on the fixed stiffness and spatial resolution. Appendix~\ref{app:inexact_newton} treats Newton convergence and the distinction between residual requirements and fixed floating-point tolerances.

\section{L-stability interval and parameter selection}\label{sec:stability}
We analyze absolute stability on the scalar linear test equation and then state its consequence for linear autonomous systems. A-stability requires well-defined stages and a complete amplification of modulus at most one on the closed left half-plane; L-stability additionally requires the amplification to vanish at infinity. Eliminating the intermediate stage gives the parameter relation from the stiff limit. Analyticity, the imaginary-axis bound, and the maximum modulus principle then determine the admissible interval. Throughout this section, $h>0$ and $C,D$ are finite real parameters.

\subsection{Model equation and elimination of the two stages}
For the scalar linear test equation $u_t=\mathcal L(u)$ with $\operatorname{Re}\lambda\le0$, set
\begin{equation}\label{eq:linearmodel}
 \mathcal{L}(u)=\lambda u.
\end{equation}
Here $u(t)$ is the scalar test solution; its numerical values retain the superscripts used for the vector method. Along this test equation,
\begin{equation}\label{eq:linearmodelrate}
 \mathcal{H}(u)=\frac{d}{dt}\mathcal{L}(u(t))=\lambda^2u.
\end{equation}
With $z=h\lambda$, the two stage equations become
\begin{equation}\label{eq:scalarstage1}
 u^{n+1/2}=u^n+\frac z4(u^n+u^{n+1/2})+\frac{z^2}{48}(u^n-u^{n+1/2}),
\end{equation}
and
\begin{equation}\label{eq:scalarstage2}
\begin{aligned}
 u^{n+1}={}&u^n+z\left[\left(\frac16+4C+\frac D2\right)u^n
            +\left(\frac23-8C+2D\right)u^{n+1/2}
            +\left(\frac16+4C-\frac{5D}{2}\right)u^{n+1}\right]\\
 &+z^2\left[C u^n+D u^{n+1/2}+\left(\frac D2-C\right)u^{n+1}\right].
\end{aligned}
\end{equation}
Solving the first equation gives
\begin{equation}\label{eq:middleamplification}
 u^{n+1/2}=g(z)u^n,
\end{equation}
where
\begin{equation}\label{eq:Gscalar}
 g(z)=\frac{1+z/4+z^2/48}{1-z/4+z^2/48}.
\end{equation}
The lower-case $g$ denotes only this intermediate amplification factor. Substituting it into the second stage gives
\begin{equation}\label{eq:fullamplification}
 u^{n+1}=R(z;C,D)u^n,
\end{equation}
with
\begin{equation}\label{eq:Rgeneral}
 R(z;C,D)=\frac{1+z\left[\left(\frac16+4C+\frac D2\right)
                  +\left(\frac23-8C+2D\right)g(z)\right]
                  +z^2[C+Dg(z)]}
 {1-\left(\frac16+4C-\frac{5D}{2}\right)z-\left(\frac D2-C\right)z^2}.
\end{equation}

\subsection{Deep-stiff limit and the relation between the parameters}
\begin{lemma}\label{lem:stifflimit}
For finite real parameters $C,D$, the complete amplification satisfies
\begin{equation}\label{eq:stifflimit}
 \lim_{\substack{|z|\to\infty\\\operatorname{Re}z\le0}}R(z;C,D)=0
\end{equation}
if and only if
\begin{equation}\label{eq:DCrelation}
 D=-C\ne0.
\end{equation}
\end{lemma}
\begin{proof}
Dividing the numerator and denominator of $g$ by $z^2/48$ gives
\[
 g(z)=\frac{1+12/z+48/z^2}{1-12/z+48/z^2},\qquad
 \frac1{1-12/z+48/z^2}=1+\frac{12}{z}+\frac{96}{z^2}+O(z^{-3}).
\]
Multiplication yields
\[
 g(z)=1+\frac{24}{z}+\frac{288}{z^2}+O(z^{-3}).
\]
Substituting this expansion into \eqref{eq:Rgeneral} gives
\[
\begin{aligned}
 \operatorname{Num}(z)
 &=(C+D)z^2+\left(\frac56-4C+\frac{53}{2}D\right)z
                  +(17-192C+336D)+O(z^{-1}),\\
 \operatorname{Den}(z)
 &=\left(C-\frac D2\right)z^2
       -\left(\frac16+4C-\frac{5D}{2}\right)z+1.
\end{aligned}
\]
In particular, $\operatorname{Num}(z)=(C+D)z^2+O(z)$; the displayed linear term is present for general parameters.

If $C-D/2\ne0$, the limit is $(C+D)/(C-D/2)$, which vanishes precisely when $C+D=0$. In that case $C-D/2=3C/2\ne0$.

The degenerate case must also be checked. If $C-D/2=0$, then $D=2C$ and $\operatorname{Den}(z)=1-(1/6-C)z$. For $C\ne0,1/6$, the numerator has leading term $3Cz^2$ and the denominator has degree one, so the ratio cannot vanish. For $C=1/6$, the denominator is exactly one and the numerator has leading term $z^2/2$, again excluding decay to zero. If $C=D=0$, direct elimination gives
\[
 R(z;0,0)=\frac{288+168z+42z^2+5z^3}{288-120z+18z^2-z^3}
 \longrightarrow-5.
\]
This case is therefore excluded as well. The rational expansions are uniform outside a sufficiently large disk, so the limit holds throughout the closed left half-plane. This proves the equivalence.
\end{proof}
Under $D=-C$, the complete function reduces to
\begin{equation}\label{eq:Rstableform}
 R(z;C)=\frac{1+z\left[\frac16+\frac72C+\left(\frac23-10C\right)g(z)\right]
                       +Cz^2[1-g(z)]}
 {1-\left(\frac16+\frac{13}{2}C\right)z+\frac32Cz^2}.
\end{equation}

\subsection{Analyticity, the imaginary axis, and the maximum modulus principle}
The remaining A-stability requirement is
\begin{equation}\label{eq:Astabilitytarget}
 \sup_{\operatorname{Re}z\le0}|R(z;C)|\le1.
\end{equation}
Let $\C_- =\{z\in\C:\operatorname{Re}z\le0\}$. It suffices to establish the following three properties:
\begin{equation}\label{eq:analyticcondition}
 R(\cdot;C)\text{ is analytic for }\operatorname{Re}z<0\text{ and continuous on }\C_-;
\end{equation}
\begin{equation}\label{eq:imaginarysup}
 \sup_{\omega\in\R}|R(i\omega;C)|\le1;
\end{equation}
\begin{equation}\label{eq:infinitylimit}
 \lim_{\substack{|z|\to\infty\\\operatorname{Re}z\le0}}R(z;C)=0.
\end{equation}
The last property follows from \cref{lem:stifflimit} under $D=-C\ne0$. We establish the first two properties below, then use the uniform decay at infinity to control the semicircular boundary in the maximum modulus argument.

Substitute \eqref{eq:Gscalar} into \eqref{eq:Rstableform} and write $R=H/K$, where
\begin{equation}\label{eq:H}
\begin{aligned}
 H(z;C)={}&\left[1+\left(\frac16+\frac72C\right)z+Cz^2\right]
                    \left(1-\frac z4+\frac{z^2}{48}\right)\\
 &+\left[\left(\frac23-10C\right)z-Cz^2\right]
                    \left(1+\frac z4+\frac{z^2}{48}\right),
\end{aligned}
\end{equation}
and
\begin{equation}\label{eq:K}
 K(z;C)=\left[1-\left(\frac16+\frac{13}{2}C\right)z+\frac32Cz^2\right]
                       \left(1-\frac z4+\frac{z^2}{48}\right).
\end{equation}

\begin{lemma}\label{lem:poles}
All zeros of the stage-denominator polynomial $K$ defined in \eqref{eq:K} lie in the open right half-plane if and only if
\begin{equation}\label{eq:polecondition}
 C\ge0.
\end{equation}
In particular, every $C>0$ makes the complete amplification analytic on a neighborhood of $\C_-$.
\end{lemma}
\begin{proof}
The second factor of \eqref{eq:K} has roots $6\pm2\sqrt3\,i$, both in the right half-plane. For $C\ne0$, the roots of the first factor of \eqref{eq:K} are
\[
 z_\pm=\frac{\frac13+13C\ \pm\sqrt{(\frac13+13C)^2-24C}}{6C}.
\]
If $C>0$ and these roots are nonreal, their real parts are $(1/3+13C)/(6C)>0$. If they are real, use both their sum and their product:
\[
 z_++z_-=\frac{2(1/6+13C/2)}{3C}>0,
 \qquad z_+z_-=\frac{2}{3C}>0.
\]
Both real roots are therefore positive as well. For $C=0$, the first factor is $1-z/6$, with the positive root $6$. For $C<0$, the real-coefficient quadratic has negative root product, so its roots are nonzero and real with opposite signs. One root lies in the open left half-plane, excluding every negative $C$ from the stated stage-denominator condition. In particular, a formula for the real part of a conjugate pair cannot be applied to these real roots. These cases prove \eqref{eq:polecondition}; the required analyticity follows from the absence of denominator zeros on $\C_-$.
\end{proof}

\begin{lemma}\label{lem:imaginaryaxis}
The imaginary-axis inequality \eqref{eq:imaginarysup} holds if and only if
\begin{equation}\label{eq:interval}
 C\in[\Ca,\Cr],\qquad
 \Ca=\frac{25-\sqrt{105}}{780},\qquad
 \Cr=\frac{25+\sqrt{105}}{780}.
\end{equation}
\end{lemma}
\begin{proof}
Direct multiplication of the factors in \eqref{eq:H}--\eqref{eq:K} gives
\begin{align*}
 H(z;C)&=1+\left(\frac7{12}-\frac{13}{2}C\right)z
            +\left(\frac7{48}-\frac{27}{8}C\right)z^2
            +\left(\frac5{288}-\frac{61}{96}C\right)z^3,\\
 K(z;C)&=1-\left(\frac5{12}+\frac{13}{2}C\right)z
            +\left(\frac1{16}+\frac{25}{8}C\right)z^2
            -\left(\frac1{288}+\frac{49}{96}C\right)z^3+\frac C{32}z^4.
\end{align*}
Putting $z=i\omega$ yields the full boundary expression
\[
 |R(i\omega;C)|^2-1
 =-\frac{3\omega^6\left(27C^2\omega^2-9360C^2+600C-8\right)}
 {(\omega^4+48\omega^2+2304)
  [81C^2\omega^4+(1521C^2-30C+1)\omega^2+36]}.
\]
The denominator is strictly positive for every real $C,\omega$. Its second factor can be written as
\[
 36\left(1-\frac32C\omega^2\right)^2
       +36\left(\frac16+\frac{13}{2}C\right)^2\omega^2.
\]
For $\omega=0$ this factor equals $36$. For $\omega\ne0$ and $1/6+13C/2\ne0$, its second summand is positive. At the remaining value $C=-1/39$, the factor is $36(1+\omega^2/26)^2>0$.

At $\omega=0$, $R(0;C)=1$. For $\omega\ne0$, positivity of the denominator makes the boundary bound equivalent to
\[
 27C^2\omega^2-9360C^2+600C-8\ge0.
\]
Necessity follows by taking $\omega\to0^+$ through nonzero values. Conversely, if the constant term is nonnegative, the nonnegative term $27C^2\omega^2$ proves the inequality at every nonzero $\omega$, while the origin already satisfies the boundary bound. Thus the necessary and sufficient condition is
\[
 -9360C^2+600C-8\ge0
 \quad\Longleftrightarrow\quad
 1170C^2-75C+1\le0.
\]
Solving this quadratic inequality gives exactly \eqref{eq:interval}. Equivalently, the unnormalized boundary difference is
\begin{equation}\label{eq:imagidentity}
 |K(i\omega;C)|^2-|H(i\omega;C)|^2
 =\frac{-1170C^2+75C-1}{3456}\omega^6+\frac{C^2}{1024}\omega^8.
\end{equation}
Both forms are retained to make the boundary calculation explicit.
\end{proof}

\begin{theorem}\label{thm:Lstable}
Within the family \eqref{eq:stage1}--\eqref{eq:stage2}, with coefficients \eqref{eq:coefficients}, the complete method is L-stable if and only if $D=-C$ and $C\in[\Ca,\Cr]$ as in \eqref{eq:interval}. Numerically $\Ca\approx0.0189141657$ and $\Cr\approx0.0451883984$.
\end{theorem}
\begin{proof}
If the method is L-stable, \cref{lem:stifflimit} gives $D=-C\ne0$. A-stability implies the imaginary-axis bound, so \cref{lem:imaginaryaxis} forces $C\in[\Ca,\Cr]$. This interval is positive, hence the stage-denominator condition in \cref{lem:poles} also holds. These facts establish necessity.

Conversely, let $C\in[\Ca,\Cr]$ and $D=-C$. Lemma~\ref{lem:poles} gives analyticity on a neighborhood of the closed left half-plane; \cref{lem:imaginaryaxis} gives modulus at most one on the imaginary axis; and \cref{lem:stifflimit} gives uniform decay to zero at infinity.

Choose $r_0$ such that $|R(z;C)|\le1$ whenever $\operatorname{Re}z\le0$ and $|z|\ge r_0$. For a fixed point $z_0$ with $\operatorname{Re}z_0<0$, choose $r>\max(r_0,|z_0|)$. On the boundary of the bounded half-disk $\{z:\operatorname{Re}z<0,\ |z|<r\}$, both the imaginary-axis segment and the semicircle have modulus at most one. The maximum modulus principle gives $|R(z_0;C)|\le1$. Continuity includes the imaginary axis. This proves A-stability on the full left half-plane, and the zero limit proves L-stability.
\end{proof}

\paragraph{Complete-step damping and the intermediate stage.}
The L-stability in \cref{thm:Lstable} concerns the completed update $R$. For any admissible $C$,
\begin{equation}\label{eq:stage_stiff_limit}
 \lim_{x\to+\infty}g(-x)=1,
 \qquad \lim_{x\to+\infty}R(-x;C)=0.
\end{equation}
For $u_t=-\kappa u$ with $u^0=1$ and fixed $h>0$, the first computed midpoint therefore tends to $1$ as $\kappa\to\infty$, whereas the exact value $e^{-\kappa h/2}$ tends to zero. The intermediate stage does not inherit the endpoint's stiff decay. This fixed-step stiff limit differs from the fixed-system limit $h\to0$ in \cref{thm:convergence}, so it is consistent with fourth-order convergence at both nodes. L-stability and the quadratic stiff damping below refer to the complete-step endpoint.

\begin{corollary}[Linear spectral interpretation]\label{prop:matrixstability}
Let $A\in\C^{N\times N}$ have spectrum in the closed left half-plane, let $h>0$ be fixed, and choose the parameters of \cref{thm:Lstable}. For $u'=Au$, both stage matrices are invertible and exact stage solves give $u^{n+1}=R(hA;C)u^n$, with every update eigenvalue of modulus at most one. If $A$ is normal, then
\begin{equation}\label{eq:normal}
 \|R(hA;C)\|_2\le1.
\end{equation}
If $A=X\Lambda X^{-1}$ is diagonalizable, then $\|R(hA;C)^n\|_2\le\kappa_2(X)$ for every integer $n\ge0$, where $\kappa_2(X)=\|X\|_2\|X^{-1}\|_2$.
\end{corollary}
\begin{proof}
Here $\mathcal H(u)=A^2u$, and the stage matrices are
\[
 M_{1/2}=I-\frac h4A+\frac{h^2}{48}A^2,\qquad
 M_1=I-\left(\frac16+\frac{13}2C\right)hA+\frac32Ch^2A^2.
\]
Their scalar denominator polynomials are the two factors in \eqref{eq:K}. By \cref{lem:poles}, neither vanishes on $h\sigma(A)$, so Schur triangularization proves invertibility. The polynomial and inverse matrices commute, and stage elimination gives
\[
 u^{n+1/2}=g(hA)u^n,\qquad u^{n+1}=R(hA;C)u^n.
\]
In a Schur form $A=UTU^*$, $R(hT;C)$ is triangular with diagonal entries $R(h\lambda_j;C)$, whose moduli are bounded by \cref{thm:Lstable}. For normal $A$, $T$ is diagonal and unitary invariance gives \eqref{eq:normal}. For diagonalizable $A$, taking norms in $R(hA;C)^n=X\operatorname{diag}(R(h\lambda_j;C)^n)X^{-1}$ gives the stated bound.
\end{proof}
This corollary covers symmetric negative semidefinite diffusion matrices. Its diagonalizable bound depends on $\kappa_2(X)$, so uniformity in stiffness or spatial resolution requires a corresponding uniform bound on this condition number. For nonnormal systems, eigenvalue inclusion alone does not imply Euclidean norm contractivity. Positivity and singular-limit preservation require additional analysis.

\subsection{Optimal choices for accuracy and deep-stiff decay}\label{sec:parameter}
The parameter is selected only after the L-stable interval has been established. The word optimal refers to a stated objective within this admissible family: minimizing the magnitude of the principal classical error coefficient, or maximizing the inverse-power order of deep-stiff attenuation.
\paragraph{Accuracy-oriented choice.}
The leading classical local error coefficient obtained from \cref{lem:defect} is
\begin{equation}\label{eq:E5}
 E_5(C)=\frac{1-75C}{2880}.
\end{equation}
The same coefficient is recovered by expanding the scalar function:
\[
 R(z;C)=e^z+E_5(C)z^5+O(z^6).
\]
The zero $C=1/75$ lies below $\Ca$. Thus $|E_5(C)|$ increases on the L-stable interval, and its minimum is attained at $\Ca$. The exact endpoint is denoted TSFO-a. 
%For the primary computational setting we retain
%\begin{equation}\label{eq:Coriginal}
% C_0=0.018924,\qquad D=-0.018924,
%\end{equation}
%and use the label TSFO-o. This value is strictly inside the admissible interval and close to $\Ca$. It is retained to reproduce the fixed-parameter experiments of the original submitted manuscript. The exact minimizer $\Ca$ is separately included as TSFO-a in the parameter comparisons. Thus $C_0$ specifies the original benchmark configuration, while the analytical minimum of $|E_5|$ remains at $\Ca$. These statements concern the stated principal coefficient and benchmark setting; they do not optimize over all methods or all finite-step errors.

\paragraph{Endpoint error constants.}
The accuracy comparison can be made explicit at both ends of the admissible interval:
\begin{equation}\label{eq:endpoint_error_constants}
 E_5(\Ca)=\frac{-73+5\sqrt{105}}{149760},\qquad
 E_5(\Cr)=\frac{-73-5\sqrt{105}}{149760}.
\end{equation}
Therefore
\[
 \frac{|E_5(\Cr)|}{|E_5(\Ca)|}
 =\frac{73+5\sqrt{105}}{73-5\sqrt{105}}\approx5.71.
\]
This recovers the smaller principal classical error at the left endpoint.

\paragraph{Damping-oriented choice.} At infinity the normalized amplification factor satisfies
\begin{equation}\label{eq:asymptotic}
 R(z;C)=\frac{5-183C}{9C}\frac1z+O(z^{-2}).
\end{equation}
The parameter choice in the following proposition is taken from the subsequent study of Huo and Su \cite{huosu}. Its cancellation argument is recalled here to connect that later damping refinement with the temporal family and its numerical comparisons.
\begin{proposition}[Unique quadratic stiff decay in the admissible family]\label{prop:quadratic}
There is a unique cancellation value in the admissible interval,
\begin{equation}\label{eq:Cq}
 \Cq=\frac5{183},\qquad D=-\frac5{183}.
\end{equation}
It is used in the subsequent PDE realization \cite{huosu}; the notation TSFO-q below denotes this setting. The maximal inverse-power decay order in this family is two.
\end{proposition}
\begin{proof}
Equation \eqref{eq:asymptotic} shows that cancellation requires $5-183C=0$, giving the unique value \eqref{eq:Cq}. Since
\[
 1170\Cq^2-75\Cq+1=-\frac{654}{3721}<0,
\]
this value lies strictly inside the admissible interval \eqref{eq:interval}. Substitution gives
\begin{equation}\label{eq:Rq}
 R(z;\Cq)=\frac{314z^2+2376z+5856}
 {5z^4-102z^3+866z^2-3480z+5856}
 =\frac{314}{5z^2}+O(z^{-3}).
\end{equation}
The numerator has nonzero quadratic coefficient. Therefore the maximal inverse-power decay order in this one-parameter L-stable subfamily is two, uniquely attained at $\Cq$. For other admissible values, the generic decay is first order.
\end{proof}

%The original choice $C_0$ keeps the leading classical error coefficient close to its minimum, while $C_q$ exchanges that smaller coefficient for quadratic stiff attenuation. The optimality in \cref{prop:quadratic} is the highest inverse-power decay order at infinity within this fixed-parameter family. At finite frequencies, the smallest complete amplification $|R|$ can occur at a different parameter. TSFO-r denotes the right endpoint only in the parameter diagnostic. 

The two objective-dependent choices are therefore
\[
\Ca=\frac{25-\sqrt{105}}{780},\qquad \Cq=\frac{5}{183},\qquad D=-C,
\]
where $\Ca$ minimizes the magnitude of the principal classical error coefficient, and $\Cq$ maximizes the inverse-power order of deep-stiff attenuation.

\section{Comparison with earlier implicit and multiderivative methods}\label{sec:relation}
\subsection{The advantages of the proposed scheme}
This paper employs the method of undetermined coefficients and direct Taylor matching to construct a fully implicit sequential two-stage fourth-order time discretization scheme based on the state-and-rate design. Its contribution lies in the complete update formulation and the characterized L-stable coefficient family, with the quadratic-decay parameter credited to the subsequent Huo--Su work \cite{huosu}. 

Compared with the related methods considered here, the proposed method offers the following three advantages in combination. Tsai et al.\ \cite{tsaicw2014} described a two-stage fourth-order semi-implicit method for the heat equation with stage order four; since the present work is a fully implicit scheme, no comparison between the two is necessary.

\noindent\textbf{(i) Improved stability at the same fourth order.}
Two-stage Gauss--Legendre and the compared classical symmetric fourth-order Hermite--Obreshkov method are A-stable but not L-stable \cite{hairerbook,mazzia2018}. The present completed update retains fourth-order accuracy and provides L-stable damping: its amplification tends to zero in the deep-stiff limit, whereas theirs tends to one. The choice $C_q=5/183$ further gives quadratic decay, $O(z^{-2})$, without adding stages.

\noindent\textbf{(ii) Higher order or fewer stages while retaining L-stability.}
Other L-stable methods already damp stiff modes; the comparison here concerns accuracy and stage count. Two-stage Radau IIA and TR--BDF2 have classical orders three and two, respectively, whereas the present method reaches order four with the same two unknown implicit stages \cite{hairerbook,hosea}. At the same fourth order, the selected L-stable SDIRK4 baseline uses five stages rather than two \cite{hairercode}. Additional trajectory-rate information thus combines L-stability with higher order at a fixed stage count, or fewer stages at a fixed order.

\noindent\textbf{(iii) Sequential implicit solution.}
The midpoint is solved first and then held fixed during the endpoint solve. Each Newton system has $N$ unknowns, avoiding the simultaneous $2N$- or $3N$-dimensional systems in the basic coupled two-stage Gauss and two-/three-stage Radau formulations. Existing single-state solvers can therefore be used successively. Spatial cells and solution components remain coupled within each stage; DIRK, TR--BDF2 and Rosenbrock methods are also sequential, and optimized collocation solvers may transform their coupled systems.

%The main costs are $\mathcal H$ and $\mathcal H_u$ evaluations, possible Hessian actions and sparse fill, and generally different factorizations for the two stages. Fewer stages do not guarantee lower runtime or smaller error at every step size. Full-step Hermite uses one unknown implicit state, while three-stage Radau IIA has the higher classical order five. The midpoint lacks deep-stiff decay; positivity and stiffness-uniform fourth order are not established. \Cref{tab:structure} summarizes the three complementary comparisons. Conservative spatial realizations are given in \cref{sec:implementation}; Sections~\ref{sec:numerics}--\ref{sec:pdenumerics} assess measured work and accuracy, including cases favoring Rosenbrock or adaptive BDF/Radau.

%\subsection{Categorized comparison and preserved properties}
\begin{table}[H]
\centering\small
\setlength{\tabcolsep}{4pt}
\renewcommand{\arraystretch}{1.12}
\caption{Complementary comparisons highlighting fourth-order accuracy, L-stability and two sequential implicit stages. Stage counts exclude the known initial state; derivative costs and comparison limits are discussed in Section~4.1.}\label{tab:structure}
\begin{tabular}{>{\raggedright\arraybackslash}p{31mm}>{\raggedright\arraybackslash}p{60mm}>{\raggedright\arraybackslash}p{58mm}}
\toprule
Comparison basis & Reference methods & Present method\\
\midrule
\textbf{I. Same fourth order} & Two-stage Gauss and symmetric Hermite: order 4, A-stable but not L-stable \cite{hairerbook,mazzia2018}. & Order 4 with L-stable endpoint damping; quadratic stiff decay at $C_q$.\\
\midrule
\textbf{II. L-stability retained} & Two-stage Radau IIA: order 3; TR--BDF2: order 2 \cite{hairerbook,hosea}. Selected L-stable SDIRK4: order 4 with five stages \cite{hairercode}. & Order 4 with two unknown implicit stages: higher order at two stages, or fewer stages at order 4, while retaining L-stability.\\
\midrule
\textbf{III. Stage solution} & Basic coupled two-stage Gauss and two-/three-stage Radau: simultaneous $2N$- or $3N$-state stage systems. & Two successive $N$-state Newton systems; no simultaneous nonlinear solve across stages.\\
\bottomrule
\end{tabular}
\end{table}

\subsection{The innovation of the proposed scheme}
The Hairer--Wanner theory \cite{hw1973} provides a general analytical framework, but does not furnish a concrete algorithm or a detailed analysis of the corresponding properties.
This subsection rewrites the already constructed method for coefficient comparison and an independent order check. %These representations do not identify it with a previously published concrete algorithm, and are not used to derive its coefficients or prove fourth-order convergence. The four Taylor equations \eqref{eq:secondconditions} have the equivalent moment form, with $\theta_i=i$ for $i\in\{0,1/2,1\}$,
%\begin{equation}\label{eq:moments}
% \sum_i a_i \theta_i^r+r\sum_i b_i \theta_i^{r-1}=\frac1{r+1},\qquad r=0,1,2,3,
%\end{equation}
%where the second sum is absent for $r=0$. 
For comparison, collect the existing states as $Y=(u^n,u^{n+1/2},u^{n+1})^T$. The same method can be written in the two-derivative form
\begin{equation}\label{eq:mdrk}
\begin{aligned}
 Y&=\e u^n+hA \mathcal{L}(Y)+h^2\widehat A \mathcal{H}(Y),\\
 u^{n+1}&=u^n+h a^T\mathcal{L}(Y)+h^2b^T\mathcal{H}(Y),
\end{aligned}
\end{equation}
where $a^T$ and $b^T$ are the last rows of
\begin{equation}\label{eq:tableaux}
 A=\begin{pmatrix}0&0&0\\c_0&c_{1/2}&0\\a_0&a_{1/2}&a_1\end{pmatrix},\qquad
 \widehat A=\begin{pmatrix}0&0&0\\d_0&d_{1/2}&0\\b_0&b_{1/2}&b_1\end{pmatrix}.
\end{equation}
The first row represents the known initial state, so there are two unknown implicit stages. The full two-parameter coefficients are retained here to identify classical members below.

%\subsection{Comparison of complete updates}\label{sec:comparison}
%Substituting the proved relation $D=-C$ into \eqref{eq:stage2} gives the complete endpoint update
%\begin{align}\label{eq:stablemethod}
% u^{n+1}={}&u^n+h\left[\left(\frac16+\frac72C\right)\mathcal{L}(u^n)+
% \left(\frac23-10C\right)\mathcal{L}(u^{n+1/2})+\left(\frac16+\frac{13}{2}C\right)\mathcal{L}(u^{n+1})\right]\notag\\
% &+Ch^2\left(\mathcal{H}(u^n)-\mathcal{H}(u^{n+1/2})-\frac32\mathcal{H}(u^{n+1})\right).
%\end{align}
%With $k=h/2$, the independently derived first stage is
%\begin{equation}\label{eq:ho}
% u^{n+1/2}=u^n+\frac k2(\mathcal{L}(u^n)+\mathcal{L}(u^{n+1/2}))+\frac{k^2}{12}(\mathcal{H}(u^n)-\mathcal{H}(u^{n+1/2})).
%\end{equation}
%This is the classical symmetric fourth-order Hermite--Obreshkov rule on a half step \cite{mazzia2018,mazzia2019,iavernaro}. This identity alone does not determine the complete update. \Cref{tab:coeffcomparison} compares endpoint weights on the common information nodes $(0,1/2,1)$; each complete method also includes its specified intermediate rule.
\begin{table}[H]\centering\small
\caption{Direct endpoint-coefficient comparison on nodes $(0,1/2,1)$. The initial state is known. Counts are active implicit states; a full-step Hermite update needs no intermediate solve. The specified intermediate rule is part of each complete method.}\label{tab:coeffcomparison}
\begin{tabular}{lccc}\toprule
Completed update & $(a_0,a_{1/2},a_1)$ & $(b_0,b_{1/2},b_1)$ & Implicit states\\\midrule
Explicit Li--Du & $(1,0,0)$ & $(1/6,1/3,0)$ & 0\\
Full-step Hermite & $(1/2,0,1/2)$ & $(1/12,0,-1/12)$ & 1\\
Two half-step Hermite & $(1/4,1/2,1/4)$ & $(1/48,0,-1/48)$ & 2\\
Present $C_q=5/183$, $D=-C_q$ & $(16/61,24/61,21/61)$ & $(5/183,-5/183,-5/122)$ & 2\\\bottomrule
\end{tabular}
\end{table}
In the unrestricted family \eqref{eq:coefficients}, $(C,D)=(1/12,0)$ gives full-step Hermite, and $(1/48,0)$ gives two half-step Hermite updates. The latter follows by adding the two half-step equations: the midpoint $\mathcal H$ terms cancel and the weights are the third row of the table. Neither pair belongs to the L-stable subfamily. Their complete responses are
\[
 R_{\rm HO}(z)=\frac{1+z/2+z^2/12}{1-z/2+z^2/12},\qquad
 R_{\rm HO,half}(z)=g(z)^2,
\]
both tending to one at infinity. By contrast, every L-stable completion tends to zero and retains nonzero midpoint weights $a_{1/2}=2/3-10C$ and $b_{1/2}=-C$.
\begin{proposition}[Distinct completed update for the compared symmetric formulas]\label{prop:distinctupdate}
No L-stable member in \cref{thm:Lstable} is identical to the classical symmetric fourth-order Hermite--Obreshkov update, the two-stage Gauss--Legendre update, or a fixed finite composition of these updates with nonzero step fractions.
\end{proposition}
\begin{proof}
Both individual classical formulas have the scalar response $R_{\rm HO}(z)$ displayed above \cite{hairerbook,mazzia2018}. It tends to one as $z\to-\infty$, also after any fixed nonzero rescaling of $z$. Every fixed finite composition has a product response with the same limiting value one. By \cref{thm:Lstable}, the proposed completed response tends to zero. Equal numerical updates would agree on every scalar linear equation, so the different limits prove the assertion. In particular, the two-half-step Hermite realization is excluded. The converse is not used: identical scalar responses alone would not establish equality of nonlinear methods.
\end{proof}

\section{Stage invertibility, Newton solution and fully discrete finite volumes}\label{sec:implementation}

\subsection{Structural preservation and stage invertibility}

\begin{proposition}\label{prop:invariants}
If a fixed row vector $\ell^T$ satisfies $\ell^T\mathcal{L}(u)=0$ for all states in the relevant neighborhood, then the exact implicit stages preserve $\ell^Tu^n$: $\ell^Tu^{n+1/2}=\ell^Tu^{n+1}=\ell^Tu^n$. An autonomous equilibrium $u_*$ with $\mathcal{L}(u_*)=0$ is also retained on the consistent stage branch.
\end{proposition}
\begin{proof}
Differentiating $\ell^T\mathcal{L}(u)=0$ with respect to the state gives $\ell^T\mathcal{L}_u(u)=0$, hence $\ell^T\mathcal{H}(u)=\ell^T\mathcal{L}_u(u)\mathcal{L}(u)=0$. Multiplication of the two stage equations by $\ell^T$ proves the invariant statement. At an equilibrium, $\mathcal{H}(u_*)=\mathcal{L}_u(u_*)\mathcal{L}(u_*)=0$, so $u^{n+1/2}=u^{n+1}=u_*$ solves both equations. Local uniqueness identifies the consistent branch. For an inexact final solve with residual $r_1$, the invariant defect in one step is exactly $\ell^Tr_1$.
\end{proof}

These shared preservation properties provide useful implementation checks.
For a linear diffusion operator $A=A^T\preceq0$, let $B=-hA\succeq0$. The two stage matrices are
\begin{equation}\label{eq:SPDmatrices}
 M_{1/2}=I+\frac14B+\frac1{48}B^2,\qquad
 M_1=I+\left(\frac16+\frac{13}2C\right)B+\frac32CB^2.
\end{equation}
For admissible $C>0$, both are symmetric positive definite and every eigenvalue is at least one. Consequently $\|M_{1/2}^{-1}\|_2,\|M_1^{-1}\|_2\le1$ for all $h>0$. This ensures stage invertibility, but not uniformly bounded condition numbers or iterative cost. For example, $B=\operatorname{diag}(0,\beta)$ gives $\kappa_2(M_{1/2})=1+\beta/4+\beta^2/48\to\infty$ as $\beta\to\infty$. Nonlinear Newton matrices additionally contain Hessian contractions.
%; see \cref{sec:implementation}.

\subsection{Sequential Newton solution}\label{sec:newton}
For $D=-C$, the residuals used by Newton iteration are
\begin{align}
 \mathcal R_{1/2}(u^{n+1/2})&=u^{n+1/2}-\frac h4\mathcal{L}(u^{n+1/2})+\frac{h^2}{48}\mathcal{H}(u^{n+1/2})
 -u^n-\frac h4\mathcal{L}(u^n)-\frac{h^2}{48}\mathcal{H}(u^n),\label{eq:resm}\\
 \mathcal R_1(u^{n+1})&=u^{n+1}-u^n
 -h\left(\frac16+\frac{13}2C\right)\mathcal{L}(u^{n+1})+\frac32Ch^2\mathcal{H}(u^{n+1})\notag\\
 &\quad-h\left(\frac16+\frac72C\right)\mathcal{L}(u^n)-h\left(\frac23-10C\right)\mathcal{L}(u^{n+1/2})\notag\\
 &\quad-Ch^2\left[\mathcal{H}(u^n)-\mathcal{H}(u^{n+1/2})\right].\label{eq:res1}
\end{align}
The midpoint Newton matrix is
\begin{equation}\label{eq:m1}
 M_{1/2}=I-\frac h4\mathcal{L}_u(u^{n+1/2})+\frac{h^2}{48}\mathcal{H}_u(u^{n+1/2}).
\end{equation}
and the endpoint matrix is
\begin{equation}\label{eq:mfinal}
 M_1=I-h\left(\frac16+\frac{13}2C\right)\mathcal{L}_u(u^{n+1})+\frac32Ch^2\mathcal{H}_u(u^{n+1}).
\end{equation}
The exact directional derivative is
\begin{equation}\label{eq:DG}
 \mathcal H_u(u)[\xi]=\mathcal L_{uu}(u)[\xi,\mathcal L(u)]
                  +\mathcal L_u(u)\bigl(\mathcal L_u(u)\xi\bigr).
\end{equation}
Thus an exact Newton implementation generally involves a Hessian action. A full third-order tensor need not be stored: a directional action or an analytically assembled contraction suffices. Automatic differentiation and matrix-free implementations are possible alternatives; their errors and costs must be included if used \cite{zeifang,chouchoulis}. The ODE implementation uses analytic derivatives and Hessian contractions. 
%For the autonomous nonlinear vector, HIRES and Robertson problems, the analytic action $\mathcal H_u(u)\xi$ is checked against a complex-step directional derivative of the separately implemented analytic $\mathcal H(u)$ at documented perturbed states and directions. The model-specific contractions and the verification procedure are supplied in the code.

%Even if the spatial matrix has a narrow sparse stencil, $\mathcal L_u^2$ may have additional nonzeros. This effect is included in the measured sparse factorizations, so fewer temporal stages do not by themselves imply a smaller spatial stencil or lower factorization cost.

The iterate $u^{n+\sigma,(\ell)}$ denotes the $\ell$th Newton approximation at the numerical stage $\sigma\in\{1/2,1\}$.
In the ordinary double-precision comparisons Newton iteration starts with $u^{n+1/2,(0)}=u^n$
and $u^{n+1,(0)}=u^{n+1/2}$. 
%In the extended-precision Robertson experiment, a
%slow predictor and the positive root of
%$3\times10^7u_2^2+10^4u_3u_2-0.04u_1=0$ provide a first-stage QSSA starting guess
%when the estimated first-step stiffness exceeds one. For smaller steps a continuous
%predictor is used. QSSA is an initialization device only: Newton still solves the
%full equations \eqref{eq:stage1} and \eqref{eq:stage2}. Each update solves the current exact Jacobian system. A backtracking reduction is used if a full update increases the residual. The stopping criterion is an infinity-norm correction bound of $3\times10^{-14}(1+\|y\|_\infty)$, with at most eighteen updates. Affine problems require one Newton update because the residual is linear. No automatic step reduction is applied to failed fixed-step runs: failures, if any, are recorded in the raw CSV output.
At either stage $k\in\{1/2,1\}$, Newton's iteration is
\begin{equation}\label{eq:newtonupdates}
\begin{aligned}
 M_k(u^{n+k,(\ell)})\delta^{(\ell)}&=-\mathcal R_k(u^{n+k,(\ell)}),\\
 u^{n+k,(\ell+1)}&=u^{n+k,(\ell)}+\alpha_\ell\delta^{(\ell)},\qquad0<\alpha_\ell\le1.
\end{aligned}
\end{equation}
The second solve begins only after the first stage has converged. A complete step is:
\begin{enumerate}[label=\textbf{\arabic*.},leftmargin=*]
\item Given $u^n$, choose admissible $C$, set $D=-C$ and evaluate $\mathcal{L}(u^n),\mathcal{H}(u^n)$ once.
\item Set $u^{n+1/2,(0)}=u^n$. For the intermediate state, evaluate $\mathcal{L},\mathcal{H},\mathcal{L}_u,\mathcal{H}_u$ at each trial vector, form $\mathcal R_{1/2},M_{1/2}$, solve \eqref{eq:newtonupdates} and apply the stopping rule. Store the converged intermediate state and its operator values.
\item Set $u^{n+1,(0)}=u^{n+1/2}$. Keep the initial and converged intermediate states and their operator values fixed. For the endpoint state, evaluate the endpoint operators and solve $\mathcal R_1=0$ by \eqref{eq:newtonupdates}.
\item Accept $u^{n+1}$ after convergence. Record residuals, evaluations, factorizations and linear solves; a failed stage is reported as a failed fixed-step run.
\end{enumerate}
%The tolerances must resolve the requested accuracy. In particular, the inexact-solve assertion of \cref{thm:convergence} uses an unscaled stage residual of order $h^5$; a fixed floating-point correction tolerance does not guarantee this asymptotically below its accuracy floor. The companion routine \code{step\_with\_stages} exposes both stage vectors and their final residuals, and is checked against the original core integrator. Appendix~\ref{app:inexact_newton} gives the residual-dependent global bound and a local Newton convergence estimate.

\subsection{WENO reconstruction and Lax--Wendroff interface evolution}\label{sec:spatialclosure}
Following Li and Du \cite{lidu}, the spatial unknowns are conservative cell averages. In each cell a fifth-order polynomial WENO reconstruction gives the point state and its gradients from the same polynomial. In one dimension the data on a common face are
\begin{equation}\label{eq:interface_data}
 \begin{aligned}
 U_{j+1/2}^{-}&=P_j(x_{j+1/2}),&
 (U_x)_{j+1/2}^{-}&=P_j'(x_{j+1/2}),\\
 U_{j+1/2}^{+}&=P_{j+1}(x_{j+1/2}),&
 (U_x)_{j+1/2}^{+}&=P_{j+1}'(x_{j+1/2}).
 \end{aligned}
\end{equation}
For example, the left state has leading terms $P_j(x_j)+\Delta x\,P_j'(x_j)/2$. The complete polynomial also retains the higher powers of the distance to the cell center. In two dimensions, transverse polynomial reconstruction supplies the strip values at three Gauss points per face, followed by the normal reconstruction and integration over the face.

The five-cell optimal quartic $P_{\rm opt}$ and three three-cell quadratics $Q_1,Q_2,Q_3$ reproduce the given averages on their stencils. The central WENO construction \cite{cravero2018} uses $\gamma=(1/2,1/6,1/6,1/6)$ and $Q_0=(P_{\rm opt}-\sum_{r=1}^3\gamma_rQ_r)/\gamma_0$. The nonlinear weights take the Z form \cite{borges2008}.
\begin{equation}\label{eq:polynomial_weno}
 P_j=\sum_{r=0}^3\frac{\alpha_r}{\sum_s\alpha_s}Q_r,\qquad
 \alpha_r=\gamma_r\left[1+\frac{(\beta_1-\beta_3)^2}{(\epsilon_W+\beta_r)^2}\right],
 \qquad \epsilon_W=10^{-12}.
\end{equation}
Here $\beta_r$ is the Jiang--Shu polynomial smoothness indicator \cite{jiangshu}. The code differentiates the same reconstruction when a discrete trajectory derivative or a transverse gradient is required.

A local Lax--Wendroff interface solve takes the two states and gradients as inputs. For the positive-speed transport and positive-state scalar Burgers fluxes used here, the local upwind state and its first time coefficient are
\begin{equation}\label{eq:interface_rate}
 U^*=U^-,\qquad
 \dot U^*=S(U^*)-\sum_{\ell=1}^{d_s}D_U f_\ell(U^*)\,(U_{x_\ell})^*.
\end{equation}
The one-sided gradients in this expression are evaluated at the selected trace. In two dimensions this includes the transverse flux derivative and the source. The face flux and its rate are $\widehat f_e=f(U^*)$ and $\widehat g_e=D_U f(U^*)\dot U^*$. For the second-order predictor used as a separate diagnostic, the midpoint face flux is
\begin{equation}\label{eq:lw2_predictor}
 F_e^{\rm LW2}=f\!\left(U_e^{*,n}+\frac h2\dot U_e^{*,n}\right).
\end{equation}
The implicit fourth-order scheme instead reconstructs the face data at each trial cell vector at $i=0,1/2,1$ and applies the weights from \eqref{eq:coefficients}. In particular,
\begin{equation}\label{eq:fv_integrated_flux}
 U_K^{n+1}=U_K^n-\frac{1}{|K|}\sum_{e\subset\partial K}|e|\,
          \mathcal F_e^n\cdot n_{K,e}+\mathcal S_K^n,\qquad
 \mathcal F_e^n=h\sum_i a_i\widehat f_e^i+h^2\sum_i b_i\widehat g_e^i .
\end{equation}
The source integral has the same weights on $\overline S$ and $\dot{\overline S}$. The midpoint balance uses $h(\widehat f_e^0+\widehat f_e^{1/2})/4+h^2(\widehat g_e^0-\widehat g_e^{1/2})/48$ and the corresponding source weights. Shared face fluxes enter adjacent cells with opposite signs. Each stage solves for all grid unknowns and includes all flux directions and the source together. The physical-interface and acoustic-control implementations use Newton--Krylov solves with complex-step Jacobian actions and a Fourier preconditioner based on the linear-weight spatial symbol. The scaled residual tolerance is $6\times10^{-14}\|U^0\|_\infty$ in the physical-interface tests and $2\times10^{-11}\|U^0\|_\infty$ in the acoustic comparison. The two stages are solved in sequence.

Write $\mathcal L_h(U)=-\operatorname{div}_h\widehat f(U)+\overline S(U)$. The two rate constructions tested below are
\begin{equation}\label{eq:closure_rates}
 \mathcal H_h^{\rm LW}=-\operatorname{div}_h\widehat g+D_U\overline S\,\mathcal L_h,
 \qquad
 \mathcal H_h^{\rm tr}=D_U\mathcal L_h\,\mathcal L_h.
\end{equation}
The classical time-order theorem applies directly to the exact discrete trajectory closure $\mathcal H_h^{\rm tr}$. The local physical interface coefficient provides a different closure on a fixed mesh. On smooth solution branches, a uniform rate difference bounded by $\delta_{\Delta x}$ enters the two stage residuals with a factor $h^2$. The local inverse and propagation bounds of Appendix~\ref{app:convergence_details} give a one-step difference of $O(h^2\delta_{\Delta x})$ and a fixed-time difference of $O(h\delta_{\Delta x})$. Consequently, the experiments distinguish fixed-mesh time refinement from joint space--time refinement of the physical-interface implementation.

\section{ODE verification: accuracy, stiff decay, and comparison}\label{sec:numerics}
The experiments follow the theoretical development. Step refinement on fixed linear and nonlinear systems, including the original Robertson configuration, first tests classical fourth order. A separate fixed-step scalar experiment then tests L-stable attenuation and inverse-power decay as stiffness grows. Autonomous reaction kinetics and an oscillator assess nonlinear convergence and cost. The two settings are TSFO-a with $C=\Ca$ and TSFO-q with $C=\Cq$, both with $D=-C$. The timing protocol is given in Section~\ref{sec:cost_protocol}.

\subsection{Original linear benchmark with separated decay rates}
We first retain the affine system
\begin{equation}\label{eq:originallinear}
 u'_1=-1000u_1+1,\qquad u'_2=-u_2+1,\qquad u(0)=(0,0)^T,
\end{equation}
at $T=10$, with $h=1/2^\ell$, $\ell=0,\ldots,9$, for both parameter choices.
Its exact solution is $u_1=(1-e^{-1000t})/1000$, $u_2=1-e^{-t}$.
To separate discretization error from floating-point subtraction, we evaluate the
exact affine numerical recurrence in 75-digit decimal arithmetic. For rate $a>0$,
the component error is $(e^{-aT}-R(-ah)^n)/a$, $n=T/h$. This is an analytically
reduced linear diagnostic, not a nonlinear solver timing experiment. The complete
sparse stage solves used for spatial tests are documented separately.
\begin{table}[t]\centering\small
\caption{Affine benchmark at $T=10$. Endpoint infinity-norm errors are evaluated from the shifted recurrence in 75-digit arithmetic.}\label{tab:linearoriginal}

% BEGIN_AUTO_linear_original
\begin{tabular}{rrrrrrr}\toprule
$h$ & TSFO-a $E_\infty$ & Order & TSFO-q $E_\infty$ & Order & GL4 $E_\infty$ & Order\\\midrule
$1.000\times10^{0}$ & $6.684\times10^{-8}$ & -- & $1.836\times10^{-7}$ & -- & $8.869\times10^{-4}$ & --\\
$5.000\times10^{-1}$ & $4.138\times10^{-9}$ & 4.014 & $1.079\times10^{-8}$ & 4.088 & $6.188\times10^{-4}$ & 0.519\\
$2.500\times10^{-1}$ & $2.580\times10^{-10}$ & 4.004 & $6.588\times10^{-10}$ & 4.034 & $1.466\times10^{-4}$ & 2.077\\
$1.250\times10^{-1}$ & $1.611\times10^{-11}$ & 4.001 & $4.075\times10^{-11}$ & 4.015 & $4.620\times10^{-7}$ & 8.310\\
$6.250\times10^{-2}$ & $1.007\times10^{-12}$ & 4.000 & $2.535\times10^{-12}$ & 4.007 & $9.624\times10^{-12}$ & 15.551\\
$3.125\times10^{-2}$ & $6.293\times10^{-14}$ & 4.000 & $1.581\times10^{-13}$ & 4.003 & $6.014\times10^{-13}$ & 4.000\\
$1.562\times10^{-2}$ & $3.933\times10^{-15}$ & 4.000 & $9.869\times10^{-15}$ & 4.002 & $3.758\times10^{-14}$ & 4.000\\
$7.812\times10^{-3}$ & $2.458\times10^{-16}$ & 4.000 & $6.165\times10^{-16}$ & 4.001 & $2.349\times10^{-15}$ & 4.000\\
$3.906\times10^{-3}$ & $1.536\times10^{-17}$ & 4.000 & $3.852\times10^{-17}$ & 4.000 & $1.468\times10^{-16}$ & 4.000\\
$1.953\times10^{-3}$ & $9.602\times10^{-19}$ & 4.000 & $2.407\times10^{-18}$ & 4.000 & $9.176\times10^{-18}$ & 4.000\\
\bottomrule\end{tabular}
% END_AUTO_linear_original
\end{table}
The coarse-step difference illustrates the effect of L-stable damping on the fast
component. The small-step portion checks fourth-order consistency on the remaining
slow component. The more extensive L-stable comparisons below test questions beyond
this known distinction between Gauss and an L-stable integrator.

\subsection{Autonomous nonlinear vector convergence}
For an independently known nonlinear trajectory, consider
\begin{equation}\label{eq:vectorproblem}
 \mathcal L(u)=\begin{pmatrix}-u_1^2\\-\kappa(u_2-u_1^2)-2u_1^3\end{pmatrix},
 \qquad u(0)=\begin{pmatrix}1\\1\end{pmatrix},\quad\kappa=10,\quad T=1.
\end{equation}
The exact solution is $u(t)=((1+t)^{-1},(1+t)^{-2})^T$. The second equation is nonlinearly coupled to the first, and the departure $u_2-u_1^2$ satisfies $(u_2-u_1^2)'=-\kappa(u_2-u_1^2)$. The analytic Jacobian and the matrix Hessian contraction are
\begin{equation}\label{eq:autonomous_vector_derivatives}
 \mathcal L_u(u)=\begin{pmatrix}-2u_1&0\\2\kappa u_1-6u_1^2&-\kappa\end{pmatrix},\qquad
 \mathcal L_{uu}(u)[\cdot,\xi]=\begin{pmatrix}-2\xi_1&0\\(2\kappa-12u_1)\xi_1&0\end{pmatrix}.
\end{equation}
All operators depend only on the state. Final-time infinity-norm errors at $8,16,32,64$ steps compare the methods under the same refinement. Maximum errors over all computed middle nodes and endpoints additionally test the two-node convergence theorem.
\begin{table}[t]\centering\small
\caption{Autonomous coupled vector problem: final-time infinity-norm errors. The last column is the observed order from 32 to 64 steps.}\label{tab:order}
% BEGIN_AUTO_order
\begin{tabular}{lrrrrr}\toprule
Method & $h=1/8$ & $h=1/16$ & $h=1/32$ & $h=1/64$ & Order\\\midrule
TSFO-a & $3.244\times10^{-7}$ & $2.052\times10^{-8}$ & $1.290\times10^{-9}$ & $8.084\times10^{-11}$ & 4.00\\
TSFO-q & $8.384\times10^{-7}$ & $5.219\times10^{-8}$ & $3.256\times10^{-9}$ & $2.034\times10^{-10}$ & 4.00\\
GL4 & $2.836\times10^{-6}$ & $1.729\times10^{-7}$ & $1.074\times10^{-8}$ & $6.705\times10^{-10}$ & 4.00\\
Radau IIA(5) & $8.608\times10^{-8}$ & $2.863\times10^{-9}$ & $9.269\times10^{-11}$ & $2.951\times10^{-12}$ & 4.97\\
SDIRK4 & $3.664\times10^{-6}$ & $2.576\times10^{-7}$ & $1.730\times10^{-8}$ & $1.125\times10^{-9}$ & 3.94\\
ROS4 & $7.223\times10^{-5}$ & $7.705\times10^{-6}$ & $6.786\times10^{-7}$ & $5.175\times10^{-8}$ & 3.71\\
HO4 & $3.188\times10^{-6}$ & $1.985\times10^{-7}$ & $1.240\times10^{-8}$ & $7.747\times10^{-10}$ & 4.00\\
Radau IIA(3) & $3.404\times10^{-5}$ & $4.497\times10^{-6}$ & $5.800\times10^{-7}$ & $7.372\times10^{-8}$ & 2.98\\
Li--Du explicit & $2.934\times10^{-5}$ & $1.347\times10^{-6}$ & $7.352\times10^{-8}$ & $4.299\times10^{-9}$ & 4.10\\
\bottomrule\end{tabular}
% END_AUTO_order
\end{table}
\begin{figure}[t]\centering\includegraphics[page=1,width=\linewidth]{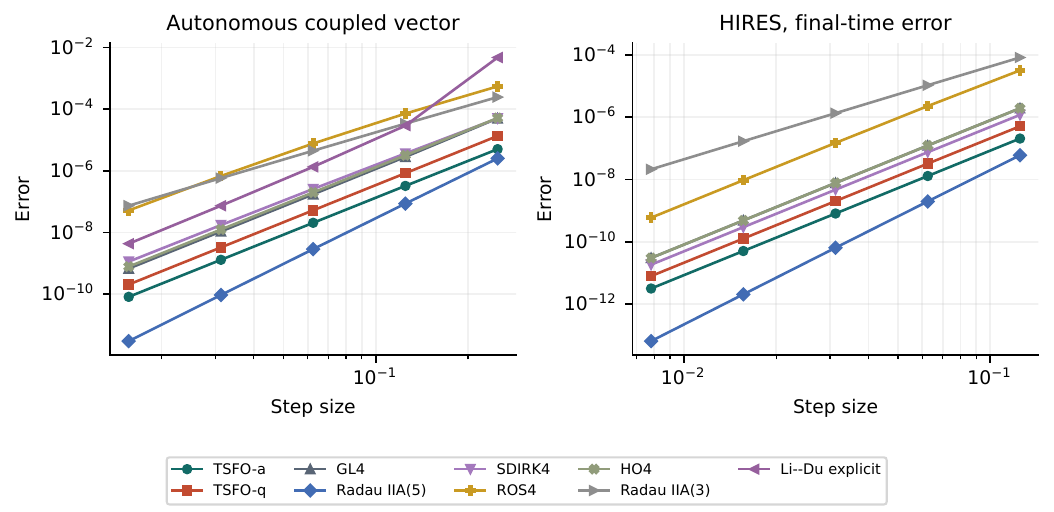}
\caption{Convergence for the autonomous nonlinear coupled vector field and HIRES. All displayed values are obtained by actual stage solves. The comparison concerns these specified systems.}\label{fig:odeaccuracy}
\end{figure}
% BEGIN_AUTO_both_nodes
The computed intermediate values are also checked against the exact trajectory, using the maximum error over all middle nodes and all endpoints of the same integration. For refinement from 32 to 64 steps, the middle-node/endpoint orders are 3.9713/4.0006 for TSFO-a and 3.9561/4.0115 for TSFO-q. At 64 steps the maximum middle-node errors are $3.307\times10^{-10}$ and $8.196\times10^{-10}$, respectively. These independently evaluated exact-orbit errors support \cref{cor:bothnodes}; the complete refinement data are retained in \code{spatial\_algorithms\_checks.json}.
% END_AUTO_both_nodes

\subsection{HIRES and Robertson chemical kinetics}\label{sec:chemistry}
The Robertson convergence study uses both parameter choices at the same step sizes.

We use the standard eight-component HIRES benchmark, not an ozone-decomposition model. With initial value $(1,0,0,0,0,0,0,0.0057)^T$, the equations are
\begin{equation}\label{eq:hires}
\begin{aligned}
 y'_1&=-1.71y_1+0.43y_2+8.32y_3+0.0007,&
 y'_2&=1.71y_1-8.75y_2,\\
 y'_3&=-10.03y_3+0.43y_4+0.035y_5,&
 y'_4&=8.32y_2+1.71y_3-1.12y_4,\\
 y'_5&=-1.745y_5+0.43y_6+0.43y_7,&
 y'_6&=-280y_6y_8+0.69y_4+1.71y_5-0.43y_6+0.69y_7,\\
 y'_7&=280y_6y_8-1.81y_7,&
 y'_8&=-280y_6y_8+1.81y_7.
\end{aligned}
\end{equation}
We specify the shortened interval $[0,1]$ for this temporal comparison; the conventional long HIRES interval is a different benchmark configuration \cite{testset,hairertest}. The final-time infinity norm is reported. This test probes a nonlinear stiff network while keeping the reference and timing experiment inexpensive enough to reproduce.

The Robertson equations are
\begin{equation}\label{eq:rober}
 y'_1=-0.04y_1+10^4y_2y_3,\quad
 y'_2=0.04y_1-10^4y_2y_3-3\times10^7y_2^2,\quad
 y'_3=3\times10^7y_2^2,
\end{equation}
with initial value $(1,0,0)^T$. The primary classical-convergence experiment retains
$T=10$ and $h=0.01/2^\ell$, $\ell=0,\ldots,9$.
At the finest step the errors are far below a reliable ordinary double-precision
subtraction of the large solution components. Therefore this experiment is recomputed
in GCC \code{\_\_float128} binary128 arithmetic (113 significand bits) \cite{quadmath}, with a Newton
correction tolerance of $10^{-30}$. It is a solution-level calculation, including the
nonlinear implicit solves, rather than a scalar error recurrence.

Mass conservation, $u_1+u_2+u_3=1$, permits exact algebraic elimination of $u_3$.
The resulting two-variable stage system is the restriction of the same full three-variable
method to its invariant plane. No QSSA is used to replace the reaction equations.
For the larger first-step stiffness, a QSSA value is used only as a Newton starting guess;
all final stage values solve the full residual. This follows the initialization strategy
in \cref{sec:implementation} without changing the integration formula.

An independent Taylor-series integration supplies the reference. Polynomial coefficient
convolutions evaluate the Robertson right-hand side; the degree-36 run with $h=2.5\times10^{-4}$
is compared with degree 40 and $h=1.25\times10^{-4}$. Both are computed in binary128.

\begin{table}[t]\centering\footnotesize
\caption{Robertson problem at $T=10$. Binary128 solutions are compared with the degree-40 Taylor reference.}\label{tab:robertsonoriginal}

% BEGIN_AUTO_robertson_original
\begin{tabular}{rrrrrrr}\toprule
$h$ & TSFO-a $E_\infty$ & Order & TSFO-q $E_\infty$ & Order & GL4 $E_\infty$ & Order\\\midrule
$1.000\times10^{-2}$ & $6.390\times10^{-10}$ & -- & $1.174\times10^{-10}$ & -- & $2.726\times10^{-9}$ & --\\
$5.000\times10^{-3}$ & $1.829\times10^{-10}$ & 1.805 & $2.269\times10^{-10}$ & -0.951 & $5.090\times10^{-10}$ & 2.421\\
$2.500\times10^{-3}$ & $5.674\times10^{-12}$ & 5.010 & $9.032\times10^{-11}$ & 1.329 & $5.041\times10^{-11}$ & 3.336\\
$1.250\times10^{-3}$ & $1.644\times10^{-12}$ & 1.787 & $5.238\times10^{-12}$ & 4.108 & $2.985\times10^{-13}$ & 7.400\\
$6.250\times10^{-4}$ & $3.308\times10^{-14}$ & 5.635 & $7.097\times10^{-14}$ & 6.206 & $1.096\times10^{-13}$ & 1.446\\
$3.125\times10^{-4}$ & $1.229\times10^{-15}$ & 4.751 & $2.920\times10^{-15}$ & 4.603 & $7.691\times10^{-15}$ & 3.833\\
$1.563\times10^{-4}$ & $6.334\times10^{-17}$ & 4.278 & $1.547\times10^{-16}$ & 4.239 & $4.934\times10^{-16}$ & 3.962\\
$7.813\times10^{-5}$ & $3.596\times10^{-18}$ & 4.138 & $8.891\times10^{-18}$ & 4.121 & $3.104\times10^{-17}$ & 3.991\\
$3.906\times10^{-5}$ & $2.140\times10^{-19}$ & 4.071 & $5.325\times10^{-19}$ & 4.062 & $1.943\times10^{-18}$ & 3.998\\
$1.953\times10^{-5}$ & $1.304\times10^{-20}$ & 4.036 & $3.257\times10^{-20}$ & 4.031 & $1.215\times10^{-19}$ & 3.999\\
\bottomrule\end{tabular}
% END_AUTO_robertson_original
\end{table}
\begin{figure}[t]\centering\includegraphics[page=2,width=\linewidth]{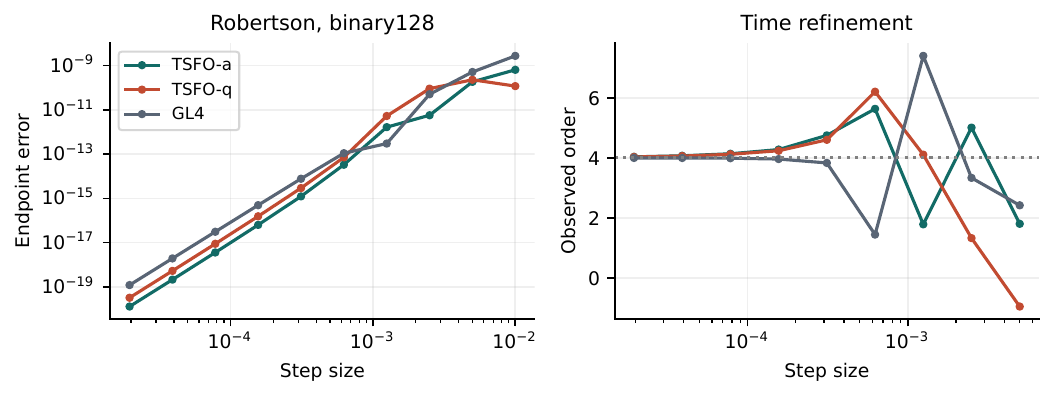}
\caption{Binary128 Robertson endpoint errors and observed refinement orders. The finest steps approach fourth order.}\label{fig:robertsonoriginal}
\end{figure}
% BEGIN_AUTO_robertson_audit
The two fine Taylor references differ by $1.107\times10^{-32}$ in the infinity norm. At the finest step, TSFO-a and TSFO-q have errors $1.304\times10^{-20}$ and $3.257\times10^{-20}$ and observed orders 4.0363 and 4.0312, respectively.
% END_AUTO_robertson_audit
The irregular orders on the coarsest steps are retained. The asymptotic trends,
not an isolated coarse-step ratio, establish the numerical fourth-order behavior.
\Cref{thm:convergence} supplies the corresponding theoretical result for this fixed
smooth polynomial system on the tested interval. This classical fourth-order result
is retained for the fixed autonomous reaction system considered here.

For the short Robertson test we use $T=0.1$ and $n=50,100,200,400,800$ for the fixed-step methods. The work--precision comparison uses the common endpoint error $E_T=\|u_{\mathrm{num}}(T)-u_{\mathrm{ref}}(T)\|_\infty$. A separate replay records the maximum componentwise error over each method's computed nodes, using the reference evaluated at those nodes. That diagnostic captures initial-layer errors that can decay by the final time. Fixed-step and adaptive runs sample their own time nodes in this diagnostic.

For the double-precision HIRES and additional early-transient Robertson comparisons, the reference uses adaptive Radau with relative tolerance $2.5\times10^{-14}$ and absolute tolerance $10^{-16}$. An independent DOP853 run with maximum step $0.002$ cross-checks the final reference value. The reference is numerical, so differences near that level are not interpreted as established convergence rates. All methods use the same analytic problem equations and starting values. The raw Robertson results also report the minimum component and the maximum mass defect along the computed trajectory; general positivity preservation is not claimed.

\begin{figure}[t]\centering\includegraphics[page=3,width=\linewidth]{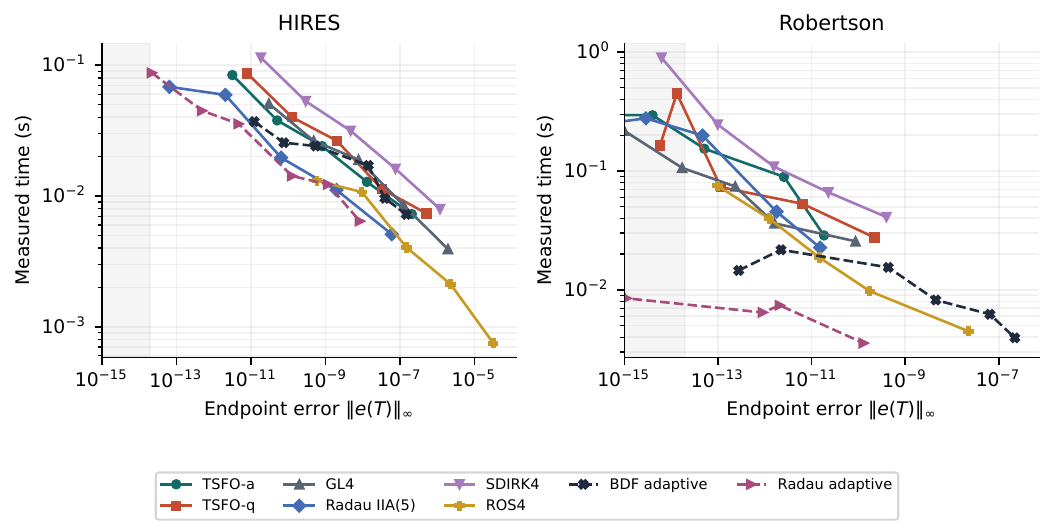}
\caption{Measured time versus endpoint infinity-norm error for HIRES at $T=1$ and Robertson at $T=0.1$. The fixed-step schemes and adaptive SciPy Radau and BDF use the same endpoint reference for each problem. Shading identifies errors near floating-point and reference uncertainty.}\label{fig:wpode}
\end{figure}

The chemical tests display why a cost comparison must accompany error tables. TSFO-a can have a smaller final-time error than GL4 at the same step, while requiring more operator work. Radau IIA(5) can be more accurate at a comparable cost in the high-accuracy regime, and ROS4 can compensate for a larger fixed-step error by its lower per-step factorization cost. \Cref{fig:wpode} reports the actual tradeoffs for these instances.

\subsection{Van der Pol oscillator}\label{sec:vanderpol}
The van der Pol oscillator \cite{hairertest} is written as
\begin{equation}\label{eq:vanderpol}
 x'=v,\qquad v'=\mu(1-x^2)v-x,\qquad \mu=100,\qquad (x(0),v(0))=(2,0).
\end{equation}
The interval $[0,1]$ examines the initial transient and classical convergence for the stated fixed problem. All stages use the analytic trajectory derivative and its exact state derivative, independently checked by complex-step directional differences. The step counts are $32,64,128,256,512,1024$, and all fixed-step implementations and Newton tolerances are the same as in the other float64 comparisons.

\begin{figure}[t]\centering\includegraphics[page=4,width=\linewidth]{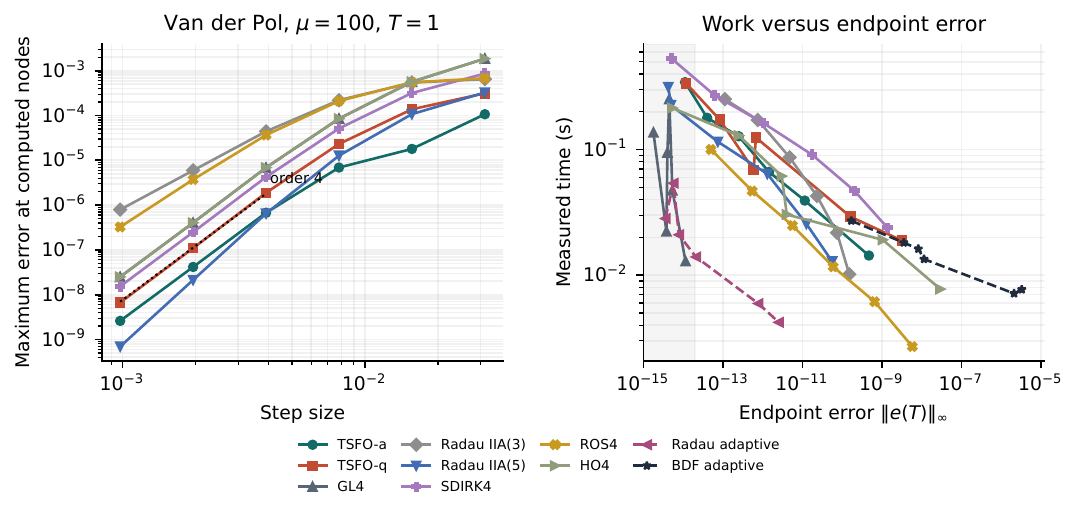}
\caption{Van der Pol at $\mu=100$ on $[0,1]$. The left panel reports the maximum error over the computed nodes under fixed-step refinement. The right panel compares measured time at a common endpoint error metric for fixed-step and adaptive methods. Shading identifies errors near floating-point and reference uncertainty.}\label{fig:vanderpol}
\end{figure}

The refinement test uses the maximum infinity-norm discrepancy over the computed nodes to measure errors in the initial transient. The work--precision comparison uses the endpoint infinity-norm error at $T=1$ for all methods. A tight Radau reference is independently checked with DOP853. Both use relative tolerance $2.5\times10^{-14}$ and absolute tolerance $10^{-16}$, with maximum steps $0.002$ and $0.0005$, respectively.
% BEGIN_AUTO_vanderpol
On 1001 common reference nodes, Radau and DOP853 differ by $3.109\times10^{-15}$ in the maximum norm. From 512 to 1024 steps, the observed orders are 3.9972 for TSFO-a and 4.0341 for TSFO-q. These rates concern the resolved initial transient of a fixed smooth system. \Cref{fig:vanderpol} uses the endpoint error at $T=1$ for the work comparison of 48 fixed-step and 12 adaptive configurations. The maximum over the computed nodes is kept as a separate transient diagnostic. The timed solves return the endpoint, and the node histories are obtained in a separate replay.
% END_AUTO_vanderpol

\subsection{Scalar L-stability and quadratic stiff decay}
Apply the actual stage solvers to $u'=-\rho u$, $u(0)=1$, with one fixed step $h=1$. Increasing $\rho$ tests stiff attenuation at fixed step, while refinement of $h$ on a fixed equation tests classical order. These are separate diagnostics. We compare the computed stage result with the independently evaluated complete amplification function. For TSFO-a and TSFO-q the latter is also evaluated in 75-digit decimal arithmetic. The explicit Li--Du formula has the fourth-degree Taylor stability polynomial and cannot remain stable as $\rho$ grows.
\begin{figure}[t]
\centering\includegraphics[page=5,width=\linewidth]{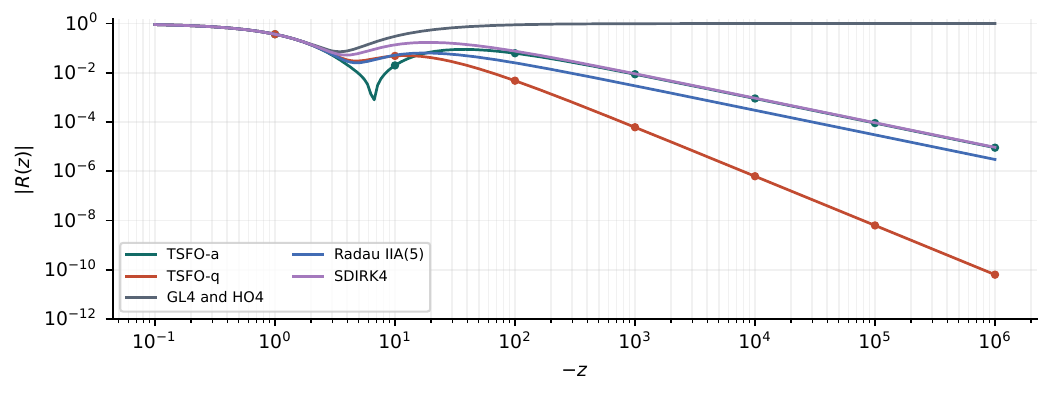}
\caption{Scalar parameter and method comparison. The plot compares the two selected parameters and implicit baselines on the negative real axis.}\label{fig:stability}
\end{figure}
\begin{table}[t]\centering\small
\caption{Magnitudes of complete-step scalar amplification functions. Values are evaluated from the implemented formulas; no time-integration timings are inferred from this table.}\label{tab:stability}
% BEGIN_AUTO_stability
\begin{tabular}{lrrrr}\toprule
Method & $z=-10$ & $z=-100$ & $z=-10^4$ & $z=-10^6$\\\midrule
Li--Du explicit & $2.910\times10^{2}$ & $4.005\times10^{6}$ & $4.165\times10^{14}$ & $4.167\times10^{22}$\\
TSFO-a & $1.987\times10^{-2}$ & $6.214\times10^{-2}$ & $9.005\times10^{-4}$ & $9.039\times10^{-6}$\\
TSFO-q & $4.833\times10^{-2}$ & $4.760\times10^{-3}$ & $6.262\times10^{-7}$ & $6.280\times10^{-11}$\\
GL4 & $3.023\times10^{-1}$ & $8.869\times10^{-1}$ & $9.988\times10^{-1}$ & $1.000\times10^{0}$\\
Radau IIA(3) & $9.589\times10^{-2}$ & $1.864\times10^{-2}$ & $1.999\times10^{-4}$ & $2.000\times10^{-6}$\\
Radau IIA(5) & $5.172\times10^{-2}$ & $2.529\times10^{-2}$ & $2.995\times10^{-4}$ & $3.000\times10^{-6}$\\
SDIRK4 & $1.366\times10^{-1}$ & $7.573\times10^{-2}$ & $9.314\times10^{-4}$ & $9.333\times10^{-6}$\\
ROS4 & $1.007\times10^{-1}$ & $2.047\times10^{-2}$ & $2.360\times10^{-4}$ & $1.740\times10^{-5}$\\
HO4 & $3.023\times10^{-1}$ & $8.869\times10^{-1}$ & $9.988\times10^{-1}$ & $1.000\times10^{0}$\\
\bottomrule\end{tabular}
% END_AUTO_stability
\end{table}

\begin{table}[t]\centering\small
\caption{Fixed-step scalar decay under increasing stiffness. Orders are logarithmic slopes with respect to $\rho$, not temporal convergence orders. Values for the sequential methods use a 75-digit evaluation independently checked against actual stage solves.}\label{tab:decayorders}
% BEGIN_AUTO_decay_orders
\begin{tabular}{rrrrr}\toprule
$\rho$ & $|R_{\mathrm a}(-\rho)|$ & Decay order & $|R_{\mathrm q}(-\rho)|$ & Decay order\\\midrule
$1.000\times10^{3}$ & $8.706\times10^{-3}$ & 0.85357 & $6.107\times10^{-5}$ & 1.89176\\
$1.000\times10^{4}$ & $9.005\times10^{-4}$ & 0.98532 & $6.262\times10^{-7}$ & 1.98908\\
$1.000\times10^{5}$ & $9.036\times10^{-5}$ & 0.99853 & $6.278\times10^{-9}$ & 1.99891\\
$1.000\times10^{6}$ & $9.039\times10^{-6}$ & 0.99985 & $6.280\times10^{-11}$ & 1.99989\\
\bottomrule\end{tabular}
% END_AUTO_decay_orders
\end{table}
\begin{figure}[t]\centering\includegraphics[page=6,width=\linewidth]{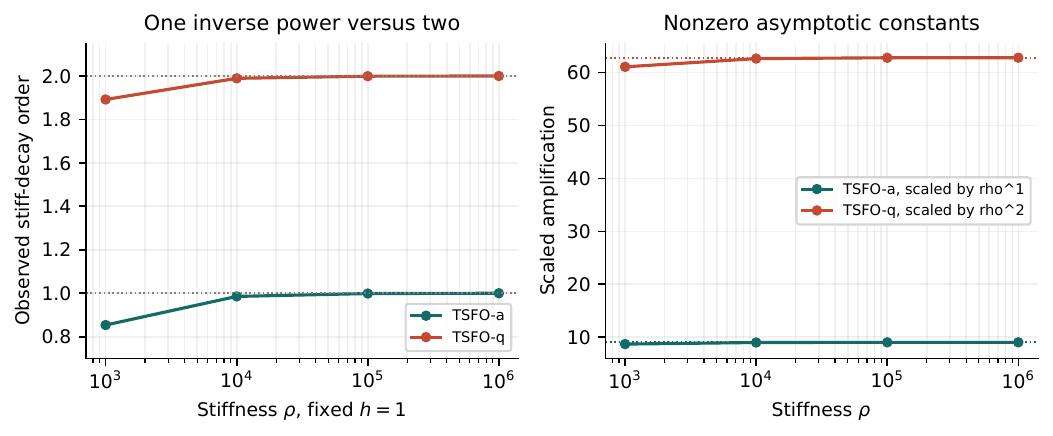}
\caption{Verification of the inverse-power decay orders. Left: observed slopes approach one for $\Ca$ and two for $\Cq$. Right: scaled amplifications approach their nonzero asymptotic constants.}\label{fig:decayorders}
\end{figure}
% BEGIN_AUTO_decay_summary
The maximum absolute discrepancy between the actual scalar stage solve and the independent TSFO amplification evaluation is $1.576\times10^{-16}$. The final decade gives decay slopes 0.99985 and 1.99989, confirming the first- and second-inverse-power predictions. This second decay order concerns $|z|\to\infty$ at fixed step; classical fourth-order time accuracy is assessed independently by step refinement.
% END_AUTO_decay_summary
The two-stage Gauss method retains a unit stiff limit, whereas the L-stable alternatives suppress the fast mode. Quadratic decay is the maximal order within the present one-parameter family. The relevant gain is enhanced decay within the same two sequential implicit solves.

\FloatBarrier
\section{PDE experiments and comparative performance}\label{sec:pdenumerics}
The tests in Section~\ref{sec:numerics} examine time convergence on fixed ODEs and attenuation under increasing stiffness at a fixed time step. We now test these properties in spatial discretizations and compare the work needed to attain prescribed errors. We use continuum cell averages to assess joint space--time convergence, same-grid references to isolate time errors, and separate measurements of the remaining fast component. Stage counts, the dimensions of the solved systems, and measured times are reported separately.

\subsection{Space--time accuracy and trajectory consistency}\label{sec:physical_fv}
\paragraph{Smooth finite-volume solutions.}
We use the WENO5 reconstruction and physical interface rate of Section~\ref{sec:spatialclosure} with both parameter choices. On periodic $[0,2\pi]$, the initial profile $p(x,0)=1+0.2\sin x$ is advanced by linear advection at speed one and by the Burgers flux $p^2/2$. The reference cell averages come from translation for advection and 24-point integration of the smooth characteristic solution for Burgers. Errors are mean absolute cell-average errors. The spatial advection study takes $T=0.1$ with $h\le 0.02\Delta x$ to keep the time error small. Joint refinement takes $T=1$ with $N/16$ steps for $N=16,32,64,128$. The separate LW2 midpoint flux uses $h\le 0.2\Delta x$.

\begin{table}[htbp]\centering\small
\caption{Physical-interface accuracy on smooth solutions. The errors are for $N=128$, with orders computed from the $64$- and $128$-cell results.}\label{tab:physical_refinement}
% BEGIN_AUTO_physical_refinement
\begin{tabular}{llrr}\toprule
Test & Method & Error & Order\\\midrule
Spatial advection & TSFO-a & $6.046\times10^{-11}$ & 4.997\\
Spatial advection & TSFO-q & $6.047\times10^{-11}$ & 4.997\\
Joint advection & TSFO-a & $4.579\times10^{-9}$ & 4.055\\
Joint advection & TSFO-q & $1.136\times10^{-8}$ & 4.015\\
Joint Burgers & TSFO-a & $5.184\times10^{-8}$ & 4.042\\
Joint Burgers & TSFO-q & $1.262\times10^{-7}$ & 3.993\\
Midpoint flux & LW2 & $2.040\times10^{-6}$ & 2.000\\
\bottomrule\end{tabular}
% END_AUTO_physical_refinement
\end{table}
The spatial advection errors approach fifth order, while joint refinement gives approximately fourth order for both equations and both TSFO choices. The LW2 midpoint flux gives second order on advection (\cref{tab:physical_refinement}).

\paragraph{Time refinement on a fixed mesh.}
Consider the Burgers transport--relaxation pair on the same periodic domain
\begin{equation}\label{eq:burgers_relaxation}
 p_t+(p^2/2)_x=0,\qquad q_t+(q^2/2)_x=-\sigma(q-p).
\end{equation}
Set $q(x,0)=p(x,0)+0.1+0.05\cos(3x)$. With $32$ cells and $\sigma=10$, integrate to $T=0.2$. For each TSFO parameter, use the same conservative residual $\mathcal L_h$ with either the local physical rate $\mathcal H_h^{\rm LW}$ or the exact discrete trajectory rate $\mathcal H_h^{\rm tr}=\mathcal L_{h,U}\mathcal L_h$. Both use the reconstruction of Section~\ref{sec:spatialclosure}. A tight DOP853 integration of $\mathcal L_h$, checked against Radau, gives the reference at each tested end-of-step node. The reported error is the maximum over those nodes, all cells, and both components.
\begin{table}[htbp]\centering\small
\caption{Maximum endpoint-node error in fixed-mesh time refinement of the coupled Burgers system. Each pair of columns uses the same first residual and changes only its rate closure.}\label{tab:closure_errors}
% BEGIN_AUTO_physical_closure
\begin{tabular}{rrrrr}\toprule
Steps & TSFO-a physical & TSFO-a trajectory & TSFO-q physical & TSFO-q trajectory\\\midrule
8 & $3.405\times10^{-6}$ & $1.585\times10^{-7}$ & $4.923\times10^{-6}$ & $5.116\times10^{-7}$\\
16 & $8.889\times10^{-7}$ & $3.483\times10^{-8}$ & $1.283\times10^{-6}$ & $1.095\times10^{-7}$\\
32 & $4.149\times10^{-7}$ & $2.348\times10^{-9}$ & $5.958\times10^{-7}$ & $6.081\times10^{-9}$\\
64 & $2.117\times10^{-7}$ & $1.596\times10^{-10}$ & $3.049\times10^{-7}$ & $4.033\times10^{-10}$\\
\bottomrule\end{tabular}
% END_AUTO_physical_closure
\end{table}
% BEGIN_AUTO_physical_closure_text
The two reference integrations agree within $8.438\times10^{-15}$. From 32 to 64 steps, the exact trajectory rate gives orders 3.879 and 3.914 for TSFO-a and TSFO-q, respectively. The corresponding orders with the physical rate are 0.971 and 0.967.
% END_AUTO_physical_closure_text
This difference is consistent with the fixed-time $O(h\delta_{\Delta x})$ contribution of the rate discrepancy derived in Section~\ref{sec:spatialclosure}. The joint refinement and fixed-mesh time refinement thus assess different accuracy properties of the physical-rate implementation.
\begin{figure}[htbp]\centering\includegraphics[page=7,width=\linewidth]{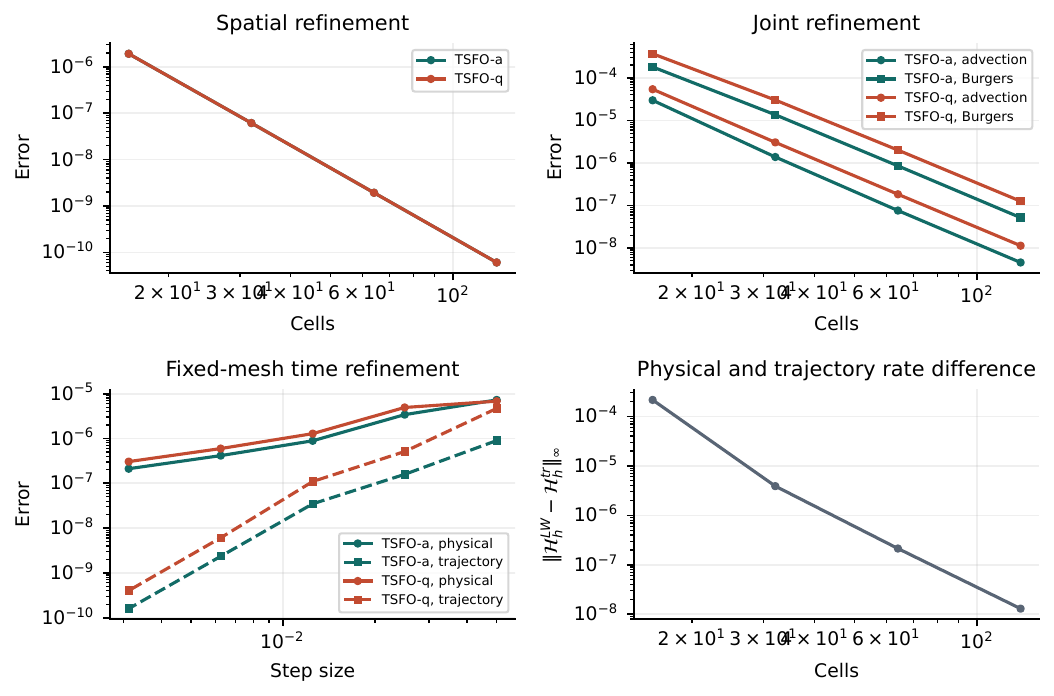}
\caption{Spatial and joint refinement with the physical interface rate. The lower panels show fixed-mesh time refinement for both rate closures and the difference between the two rates under spatial refinement.}\label{fig:physical_refinement}
\end{figure}
\FloatBarrier

\paragraph{Longer integrations and fixed-grid diffusion.}
For comparison with the shorter ODE tests in Section~\ref{sec:numerics}, we integrate HIRES to $T=1,10,321.8122$ and van der Pol at $\mu=100$ to $T=100$. The latter uses binary128 arithmetic with $h=0.01/2^\ell$, $\ell=0,\ldots,4$. We also test the heat and biharmonic equations on $(0,\pi)^d$, $d=1,2$, using $12$ interior points per direction and initial data $\prod_{j=1}^d\sin(x_j)$. The discrete operators are the centered Dirichlet Laplacian $A_\Delta$ and $-A_\Delta^2$. The heat problem has homogeneous Dirichlet data and $T=1$. The biharmonic problem has $u=\Delta u=0$ at the boundary and $T=0.1$. The stage solves in the discrete sine basis are checked against direct sparse solves, with the exact discrete exponential as the time reference.
\begin{table}[htbp]\centering\small
\caption{Endpoint errors and observed orders on the finest two time steps. Spatial grids are held fixed in the diffusion tests.}\label{tab:extended_accuracy}
% BEGIN_AUTO_extended_accuracy
\begin{tabular}{lrrrrr}\toprule
Test & Finest $h$ & TSFO-a error & Order & TSFO-q error & Order\\\midrule
HIRES, $T=1$ & $1.562\times10^{-2}$ & $5.048\times10^{-11}$ & 4.000 & $1.265\times10^{-10}$ & 3.999\\
HIRES, $T=10$ & $1.562\times10^{-2}$ & $4.149\times10^{-12}$ & 3.997 & $1.045\times10^{-11}$ & 4.002\\
HIRES, $T=321.812$ & $6.246\times10^{-2}$ & $3.336\times10^{-10}$ & 4.004 & $8.494\times10^{-10}$ & 4.005\\
Van der Pol, $T=100$ & $6.250\times10^{-4}$ & $6.630\times10^{-9}$ & 3.998 & $1.639\times10^{-8}$ & 3.976\\
Heat, 1D & $6.250\times10^{-2}$ & $7.943\times10^{-10}$ & 4.000 & $2.000\times10^{-9}$ & 4.007\\
Heat, 2D & $6.250\times10^{-2}$ & $9.329\times10^{-9}$ & 4.001 & $2.359\times10^{-8}$ & 4.015\\
Biharmonic, 1D & $6.250\times10^{-3}$ & $1.899\times10^{-14}$ & 4.000 & $4.763\times10^{-14}$ & 4.001\\
Biharmonic, 2D & $6.250\times10^{-3}$ & $1.434\times10^{-11}$ & 4.000 & $3.602\times10^{-11}$ & 4.003\\
\bottomrule\end{tabular}
% END_AUTO_extended_accuracy
\end{table}
% BEGIN_AUTO_extended_text
The largest independent endpoint reference gap for the nonlinear systems is $1.887\times10^{-15}$. For the diffusion problems the four-step sparse and sine-basis solutions agree within $5.117\times10^{-14}$. The errors in the diffusion rows isolate the time discretization.
% END_AUTO_extended_text
The orders in \cref{tab:extended_accuracy} are close to four for both parameters.
\FloatBarrier

\subsection{Fast-mode damping and the two parameter choices}
\subsubsection{Accuracy and damping objectives}
The two settings have the same formal order and use the same two sequential implicit stages. The comparison in \cref{tab:parameter_comparison} reports the magnitude of the principal error coefficient, the HIRES endpoint error at $T=1$ with $h=1/64$, and three measured full-step gains. The gains are computed by solving the stages for a real two-component representation of $u'=zu$ and checked against the scalar amplification function.
\begin{table}[htbp]\centering\small
\caption{Principal time-error coefficient, HIRES endpoint error, and measured full-step gains for the two parameter choices.}\label{tab:parameter_comparison}
% BEGIN_AUTO_parameter_comparison
\begin{tabular}{lrrrrrr}\toprule
Setting & $C$ & $|E_5(C)|$ & HIRES error & $|R(-10)|$ & $|R(-1000)|$ & $|R(5i)|$\\\midrule
TSFO-a & 0.01891417 & $1.453\times10^{-4}$ & $5.048\times10^{-11}$ & $1.987\times10^{-2}$ & $8.706\times10^{-3}$ & 0.982390\\
TSFO-q & 0.02732240 & $3.643\times10^{-4}$ & $1.265\times10^{-10}$ & $4.833\times10^{-2}$ & $6.107\times10^{-5}$ & 0.892613\\
\bottomrule\end{tabular}
% END_AUTO_parameter_comparison
\end{table}
% BEGIN_AUTO_parameter_text
At the same step size the HIRES error for TSFO-a is 0.40 times that for TSFO-q. At $z=-1000$, TSFO-q reduces the remaining amplitude by a factor of 142.6 relative to TSFO-a. At $z=-10$, TSFO-a leaves a smaller amplitude, and at $z=5i$ it retains a larger oscillatory amplitude.
% END_AUTO_parameter_text

\subsubsection{Damping with the physical interface rate}
For \eqref{eq:burgers_relaxation} we take one step of length $h=0.1$ and vary $\sigma$ over $10^2,10^3,10^4,10^5$. Fluxes cancel in the periodic sum, so the mean of $q-p$ measures the stiff relaxation directly. We report its numerical multiplier $\eta=|\overline{q-p}(h)/\overline{q-p}(0)|$. At $\sigma h=1000$, the full endpoint error $E_T$ is the maximum absolute difference from a separate same-mesh Radau reference. In this test TSFO-q reduces both the fast mean and the full endpoint error relative to Gauss, two-stage Radau IIA, and SDIRK4 (\cref{tab:physical_stiff}).
\begin{table}[htbp]\centering\small
\caption{One-step Burgers transport--relaxation at $\sigma h=1000$. The endpoint error, fast mean multiplier, and Newton work are separate diagnostics.}\label{tab:physical_stiff}
% BEGIN_AUTO_physical_stiff
\begin{tabular}{lrrrr}\toprule
Method & $E_T$ & $\eta$ & Newton updates & Max. system\\\midrule
TSFO-a & $1.298\times10^{-3}$ & $8.706\times10^{-3}$ & 5 & 64\\
TSFO-q & $9.107\times10^{-6}$ & $6.107\times10^{-5}$ & 5 & 64\\
GL4 & $1.473\times10^{-1}$ & $9.881\times10^{-1}$ & 3 & 128\\
Radau IIA(3) & $2.971\times10^{-4}$ & $1.986\times10^{-3}$ & 3 & 128\\
SDIRK4 & $1.362\times10^{-3}$ & $9.138\times10^{-3}$ & 10 & 64\\
\bottomrule\end{tabular}
% END_AUTO_physical_stiff
\end{table}

The two-dimensional physical-rate test uses the unsplit system
\begin{equation}\label{eq:fv2dtest}
 p_t+ap_x+bp_y=0,\qquad q_t+aq_x+bq_y=-\sigma(q-p)
\end{equation}
with periodic boundaries on $[0,2\pi]^2$, $(a,b)=(1,1/2)$, $\sigma=10^3$, and a $12\times12$ grid. The initial cell averages come from $p=1+0.1\sin(x+y)$ and $q-p=0.1+0.05\cos(2x-y)$. The reference at $T=0.1$ uses translation and exponential decay of the mismatch. The one-step fast gain is the ratio of the maximum absolute cell-average mismatch at $T$ to its initial value.
% BEGIN_AUTO_physical_stiff_text
The last stiffness decade gives fast mean decay exponents 0.985 and 1.989 for TSFO-a and TSFO-q. The measured one-step two-dimensional mismatch ratios are $6.215\times10^{-2}$, $4.761\times10^{-3}$, and $8.869\times10^{-1}$ for TSFO-a, TSFO-q, and Gauss, respectively.
% END_AUTO_physical_stiff_text
\begin{figure}[htbp]\centering\includegraphics[page=8,width=\linewidth]{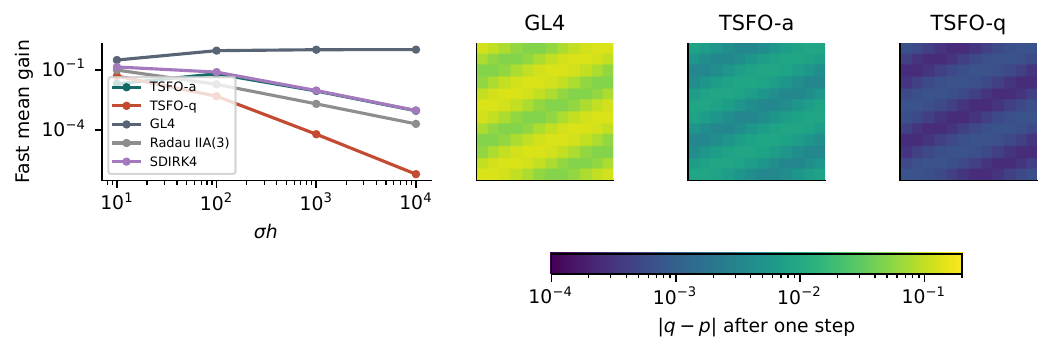}
\caption{Physical-interface tests. The left panel shows the fast mean multiplier in the Burgers stiffness sweep. The remaining panels show two-dimensional fast mismatch fields after one step on a common logarithmic scale.}\label{fig:physical_stiffness}
\end{figure}
\FloatBarrier

\subsubsection{One-dimensional conservative trajectory-rate test}\label{sec:fvtransport}
To isolate the temporal update with an exact discrete trajectory rate, consider the periodic one-dimensional system on $[0,2\pi]$
\begin{equation}\label{eq:fvpde}
 p_t+a p_x=0,\qquad q_t+a q_x=-\sigma(q-p),\qquad a=1,\quad\sigma=10^4.
\end{equation}
On $n=128$ cells of width $\Delta x=2\pi/n$, use the fourth-order centered face flux
\begin{equation}\label{eq:fvface}
 \widehat f_{j+1/2}(U)=\frac a{12}(-U_{j-1}+7U_j+7U_{j+1}-U_{j+2}).
\end{equation}
For $\widehat g_{j+1/2}$ the same face formula acts on $\dot U=\mathcal L_h(U)$, and the source rate is $-\sigma(\dot q-\dot p)$. If $D_h$ denotes the centered difference induced by \eqref{eq:fvface}, the independently assembled operators satisfy
\begin{equation}\label{eq:fvblock}
 A_h=\begin{pmatrix}-aD_h&0\\\sigma I&-aD_h-\sigma I\end{pmatrix},\qquad
 \mathcal L_h(U)=A_h U,\quad \mathcal H_h(U)=A_h^2 U.
\end{equation}
The stages are solved for all $256$ unknowns. The constant Newton matrices are formed at setup, and the rate evaluations during time integration use the shared face operations.

The initial data are the cell averages of $p=1+0.2\sin x$ and $q-p=0.2\cos(12x)$. The exact discrete reference uses the cell-average sinc factors, frequency
\[
 \omega_h(k)=\frac{a}{6\Delta x}\{8\sin(k\Delta x)-\sin(2k\Delta x)\},
\]
and the fast multiplier $e^{-\sigma t}$. The reported errors at $T=0.1$ are $\sqrt{\Delta x}\|U-U_e\|_2$, and the fast amplitude is $\sqrt{\Delta x/\pi}\|q-p\|_2$. The reference fast amplitude is recorded on a logarithmic scale below the floating-point range. Fixed-step implicit runs use $1,2,4,8,16,32$ steps. The explicit Li--Du method also uses $512$ and $1024$ steps. We check the maximum gain over the full discrete spectrum before including explicit runs in the work--precision comparison.
\begin{figure}[htbp]\centering\includegraphics[page=9,width=\linewidth]{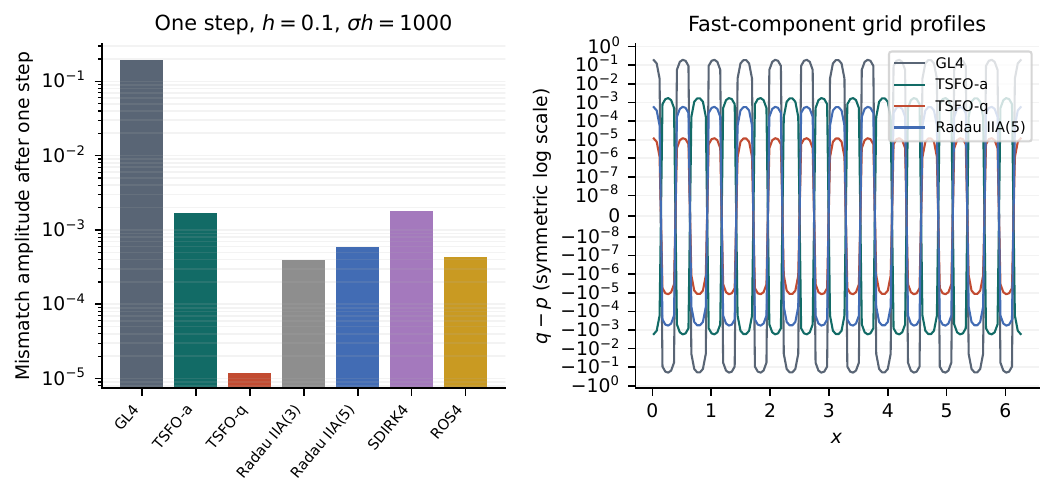}
\caption{Fast-mode attenuation in the one-dimensional finite-volume system at $h=0.1$. The left panel gives mismatch amplitudes, and the right panel gives grid profiles on a symmetric logarithmic scale.}\label{fig:fvfast}
\end{figure}
\begin{figure}[htbp]\centering\includegraphics[page=10,width=\linewidth]{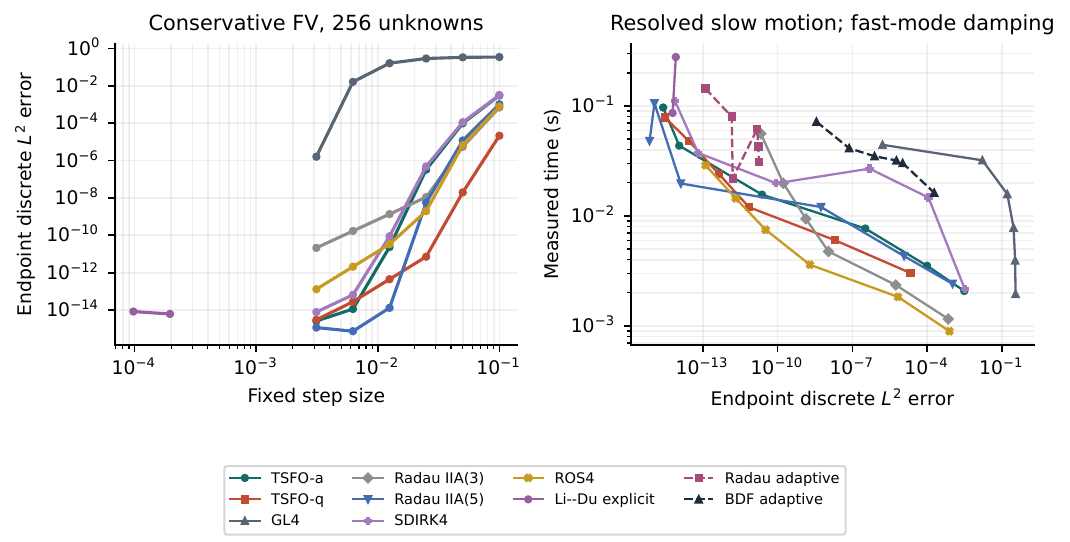}
\caption{Time refinement and work--precision for one-dimensional transport--relaxation. The fixed-step explicit curve uses only full-spectrum stable runs. Coarse explicit results are retained in \cref{tab:fvcost} and the raw data.}\label{fig:fvprecision}
\end{figure}
\begin{table}[htbp]\centering\small
\caption{One-dimensional transport--relaxation at one coarse fixed step, with a stable explicit run and selected adaptive runs attaining endpoint error $10^{-8}$. Times are medians of three runs including face, rate, and linear solve costs.}\label{tab:fvcost}
% BEGIN_AUTO_fv_cost
\begin{tabular}{lrrrrr}\toprule
Method & Steps & Error & Fast amplitude & Time (s) & Max. gain\\\midrule
TSFO-a & 1 & $3.042\times10^{-3}$ & $1.716\times10^{-3}$ & 0.0021 & $1.000\times10^{0}$\\
TSFO-q & 1 & $2.134\times10^{-5}$ & $1.204\times10^{-5}$ & 0.0030 & $1.000\times10^{0}$\\
GL4 & 1 & $3.452\times10^{-1}$ & $1.948\times10^{-1}$ & 0.0020 & $1.000\times10^{0}$\\
Radau IIA(3) & 1 & $6.939\times10^{-4}$ & $3.915\times10^{-4}$ & 0.0012 & $1.000\times10^{0}$\\
Radau IIA(5) & 1 & $1.030\times10^{-3}$ & $5.814\times10^{-4}$ & 0.0024 & $1.000\times10^{0}$\\
SDIRK4 & 1 & $3.193\times10^{-3}$ & $1.801\times10^{-3}$ & 0.0022 & $1.000\times10^{0}$\\
ROS4 & 1 & $7.715\times10^{-4}$ & $4.353\times10^{-4}$ & 0.0009 & $1.000\times10^{0}$\\
HO4 & 1 & $3.452\times10^{-1}$ & $1.948\times10^{-1}$ & 0.0016 & $1.000\times10^{0}$\\
Li--Du explicit & 1 & $1.450\times10^{10}$ & $8.181\times10^{9}$ & 0.0003 & $4.150\times10^{10}$\\
Li--Du explicit & 512 & $6.162\times10^{-15}$ & $0$ & 0.0864 & $1.000\times10^{0}$\\
\midrule
Radau adaptive & 12 & $1.592\times10^{-12}$ & $5.443\times10^{-14}$ & 0.0219 & --\\
BDF adaptive & 175 & $3.607\times10^{-9}$ & $1.093\times10^{-16}$ & 0.0719 & --\\
\bottomrule\end{tabular}
% END_AUTO_fv_cost
\end{table}
% BEGIN_AUTO_fv_summary
At the coarse step $h=0.1$, the measured mismatch amplitudes are $1.716\times10^{-3}$ for TSFO-a, $1.204\times10^{-5}$ for TSFO-q, and $1.948\times10^{-1}$ for Gauss. The explicit Li--Du step has maximum spectral amplification $4.150\times10^{10}$ and is unstable at this step size. The largest measured invariant defect in the implicit runs is $1.395\times10^{-15}$. Face-based $\mathcal H_h$ and independently assembled $A_h^2U$ agree to relative discrepancy $1.980\times10^{-16}$. For the sampled $10^{-8}$ endpoint-error target, among spectrally stable fixed-step configurations and completed adaptive runs the least measured costs are TSFO-a (8 steps, 0.0156 s), TSFO-q (4 steps, 0.0120 s), Li--Du explicit (512 steps, 0.0864 s), Radau adaptive (12 steps, 0.0219 s), ROS4 (4 steps, 0.0036 s). These times compare the accuracy of the endpoint after the relaxation layer has decayed.
% END_AUTO_fv_summary
\FloatBarrier

\subsubsection{Unsplit two-dimensional trajectory-rate solve}\label{sec:fv2dtest}
We also solve \eqref{eq:fv2dtest} with the shared centered face fluxes in both directions on a $16\times12$ periodic grid, using $(a,b)=(1,1/2)$, $\sigma=10^4$, and $T=0.1$. The full rate applies the linear residual twice, including both spatial directions and the source in each application. For cell centers $(x_i,y_j)$, let $\eta_{kl}=\operatorname{sinc}(k\Delta x/2)\operatorname{sinc}(l\Delta y/2)$ with $\operatorname{sinc}(s)=\sin(s)/s$. The initial averages are
\begin{align*}
 p_{ij}(0)&=1+0.1\eta_{11}\sin(x_i+y_j)+0.05\eta_{2,-3}\cos(2x_i-3y_j),\\
 q_{ij}(0)&=p_{ij}(0)+0.2\eta_{5,-4}\cos(5x_i-4y_j).
\end{align*}
The exact reference advances each mode with frequency
\[
 \omega_{kl}=a\frac{8\sin(k\Delta x)-\sin(2k\Delta x)}{6\Delta x}
             +b\frac{8\sin(l\Delta y)-\sin(2l\Delta y)}{6\Delta y},
\]
with an additional factor $e^{-\sigma t}$ for $q-p$. Each time integrator solves the full $384$-unknown spatial problem. The endpoint error is the discrete $L^2$ norm over both components, and the fast gain is $\|q-p\|_2/\|q(0)-p(0)\|_2$. The eight fixed-step methods use $1,2,4,8,16,32,64$ steps. The companion data include the face/matrix identity checks, derivative action checks, and stage residuals.
\begin{figure}[htbp]\centering\includegraphics[page=11,width=\linewidth]{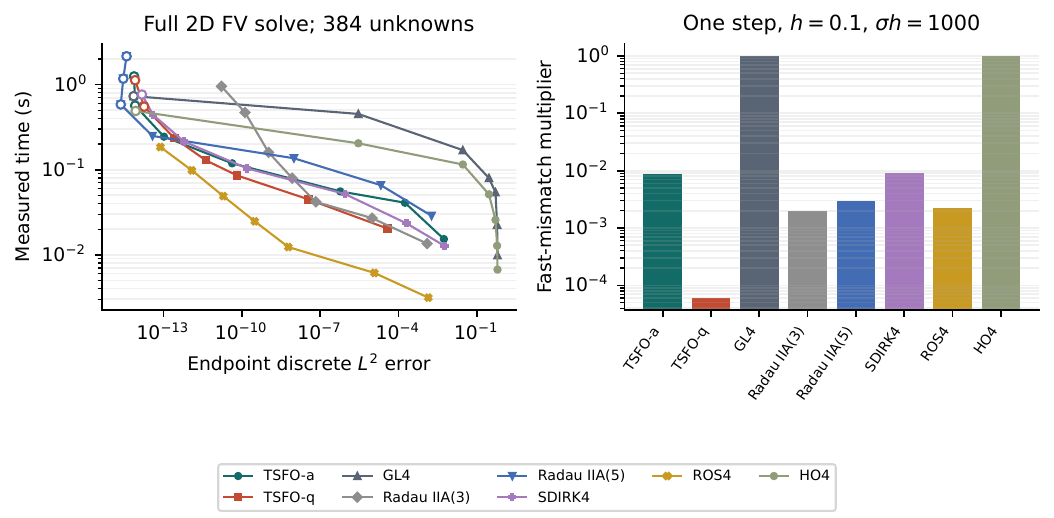}
\caption{Work--precision and one-step fast-mode attenuation for the unsplit two-dimensional conservative problem. The same discrete residual and reference are used for all methods.}\label{fig:spatialalgorithm}
\end{figure}
% BEGIN_AUTO_spatial_algorithm
All 56 configurations completed successfully. After one step, the measured fast multipliers are $8.706\times10^{-3}$ for TSFO-a, $6.107\times10^{-5}$ for TSFO-q, and $9.881\times10^{-1}$ for Gauss. The largest conserved-$p$ mass defect is $2.580\times10^{-15}$. At the endpoint error target $10^{-8}$, the smallest sampled TSFO-q time is $8.531\times10^{-2}$ seconds with 4 steps; ROS4 takes $1.232\times10^{-2}$ seconds with 4 steps. Hollow symbols mark results at the floating-point scale, below about $3\times10^{-14}$.
% END_AUTO_spatial_algorithm
\FloatBarrier

\subsubsection{High-frequency damping in a sparse heat solve}
\begin{figure}[htbp]\centering\includegraphics[page=12,width=\linewidth]{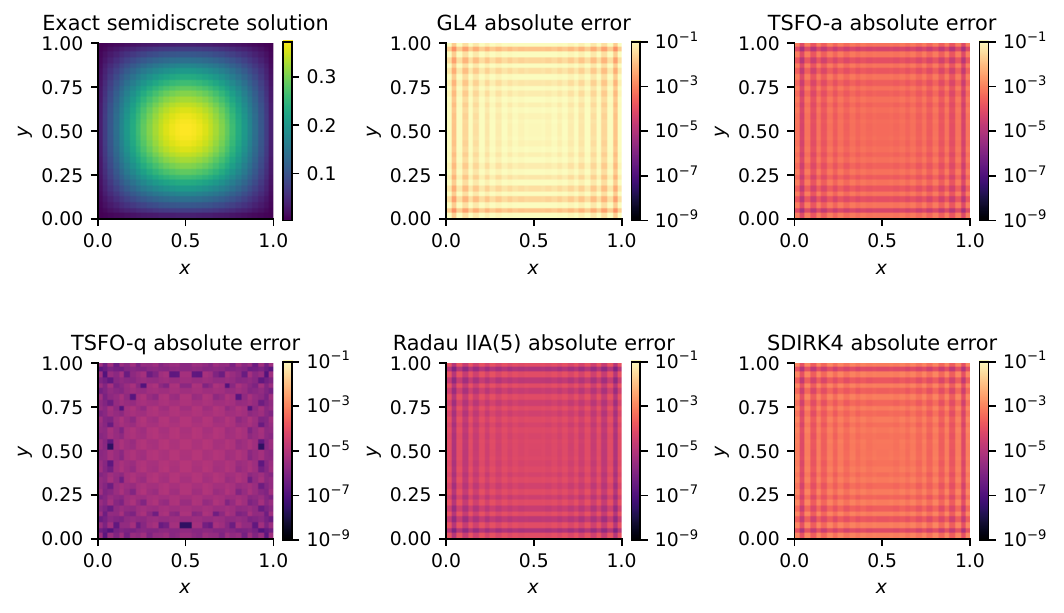}
\caption{Sparse heat equation with two steps to $T=0.05$. The first panel shows the exact semidiscrete solution. The other panels show absolute errors on a common logarithmic scale.}\label{fig:heat}
\end{figure}

On the interior $32\times32$ grid of $(0,1)^2$ with homogeneous Dirichlet data, take $U'=A_hU$ with the five-point Laplacian, spacing $\Delta x=1/33$, and initial grid data sampled from
\begin{equation}\label{eq:heatinitial}
 \sin(\pi x)\sin(\pi y)+0.2\sin(16\pi x)\sin(17\pi y).
\end{equation}
The exact semidiscrete reference uses the eigenvalues
\begin{equation}\label{eq:discreteeig}
 \lambda_{pq,h}=-\frac4{\Delta x^2}
 \left[\sin^2\!\left(\frac{p\pi\Delta x}{2}\right)
      +\sin^2\!\left(\frac{q\pi\Delta x}{2}\right)\right].
\end{equation}
Time integration uses the full sparse stage systems with $1,2,4,8$ steps to $T=0.05$. The two-step fields in \cref{fig:heat} show the high-frequency component left by Gauss and its attenuation by the L-stable integrators. The remaining field error also includes the time-discretization error of the low-frequency component.
\FloatBarrier

\subsection{Implicit stage counts and system size}
\subsubsection{Sequential and stacked solves of the same equations}
The midpoint equation is independent of the endpoint unknown, which gives a zero upper-right block in the stacked Jacobian. To check that sequential solution preserves the stage values, we solve the same pair of equations on HIRES at $h=1/4$ both sequentially and as a single stacked Newton problem. \Cref{tab:stage_agreement} shows agreement to roundoff while the sequential formulation solves systems of dimension $8$ instead of $16$.
\begin{table}[htbp]\centering\small
\caption{Sequential and stacked Newton solves for one HIRES step. The stage difference is the infinity norm over the two stage vectors.}\label{tab:stage_agreement}
% BEGIN_AUTO_sequential_table
\begin{tabular}{lrrrr}\toprule
Setting & Sequential dimension & Stacked dimension & Stage difference & Residual\\\midrule
TSFO-a & 8 & 16 & $1.110\times10^{-16}$ & $1.388\times10^{-17}$\\
TSFO-q & 8 & 16 & $1.110\times10^{-16}$ & $1.110\times10^{-16}$\\
\bottomrule\end{tabular}
% END_AUTO_sequential_table
\end{table}

\subsubsection{Total implicit stages at prescribed error}
For HIRES, all methods use the same initial value, final time $T=1$, and endpoint reference. The reported count is the number of unknown implicit stages per step times the fewest sampled dyadic steps attaining the target. The known initial state is excluded. TR--BDF2 uses a trapezoidal substep with $\gamma=2-\sqrt{2}$ followed by the BDF2 completion \cite{hosea}. The highest system dimension in \cref{tab:hires_matched} distinguishes the coupled Radau stages from the sequential solves.
\begin{table}[H]\centering\small
\caption{Total implicit stage values required to meet each sampled HIRES error target. Parentheses give the corresponding step counts.}\label{tab:hires_matched}
% BEGIN_AUTO_hires_matched
\begin{tabular}{lrrrrr}\toprule
Method & Stages/step & Max. system & $10^{-6}$ & $10^{-8}$ & $10^{-10}$\\\midrule
TSFO-a & 2 & 8 & 16 (8) & 64 (32) & 128 (64)\\
TSFO-q & 2 & 8 & 16 (8) & 64 (32) & 256 (128)\\
Radau IIA(3) & 2 & 16 & 128 (64) & 512 (256) & 2048 (1024)\\
Radau IIA(5) & 3 & 24 & 24 (8) & 48 (16) & 96 (32)\\
SDIRK4 & 5 & 8 & 80 (16) & 160 (32) & 640 (128)\\
TR-BDF2 & 2 & 8 & 512 (256) & 8192 (4096) & 65536 (32768)\\
\bottomrule\end{tabular}
% END_AUTO_hires_matched
\end{table}
% BEGIN_AUTO_hires_matched_text
The minimum sampled count at target $10^{-6}$ is 16 stage values for TSFO-a and TSFO-q. The minimum sampled count at target $10^{-8}$ is 48 stage values for Radau IIA(5). The minimum sampled count at target $10^{-10}$ is 96 stage values for Radau IIA(5).
% END_AUTO_hires_matched_text

For a two-dimensional comparison with nonlinear spatial weights, solve the periodic relaxing acoustic system
\[
 p_t+\nabla\cdot v=0,\qquad v_t+\nabla(p+0.6q)=0,\qquad
 q_t+0.6\nabla\cdot v=-4000q
\]
on $[0,2\pi]^2$ with a $32\times32$ grid and final time $T=0.5$. The initial slow and fast Fourier modes have wave numbers $(1,1)$ and $(2,3)$ and amplitudes $0.002$ and $0.0015$. All methods use the same characteristic WENO-Z5 residual with the acoustic interface flux. TSFO uses the exact rate along this discrete residual, including the differentiation of the WENO weights. The relative RMS endpoint error is normalized by the RMS of the initial data. The reference uses 256 Radau IIA(5) steps, with a check at 128 steps. We compare $1,2,4,8,16$ equal steps at a common error threshold of $10^{-6}$.
\begin{table}[H]\centering\small
\caption{Two-dimensional acoustic comparison at relative RMS temporal error $\le 10^{-6}$. Each row uses the fewest sampled steps meeting the target. Times are medians of three runs.}\label{tab:matched_pde}
% BEGIN_AUTO_weno_budget
\begin{tabular}{lrrrrrr}\toprule
Method & Steps & Stages & Error & Linear solves & Max. system & Seconds\\\midrule
TSFO-a & 4 & 8 & $3.168\times10^{-7}$ & 10 & 4096 & 0.427\\
TSFO-q & 4 & 8 & $4.846\times10^{-7}$ & 10 & 4096 & 0.428\\
Radau IIA(3) & 16 & 32 & $6.780\times10^{-7}$ & 17 & 8192 & 0.458\\
Radau IIA(5) & 4 & 12 & $1.356\times10^{-8}$ & 5 & 12288 & 0.215\\
SDIRK4 & 4 & 20 & $4.711\times10^{-7}$ & 21 & 4096 & 0.295\\
GL4 & -- & -- & -- & -- & -- & --\\
\bottomrule\end{tabular}
% END_AUTO_weno_budget
\end{table}
% BEGIN_AUTO_weno_budget_text
The relative reference refinement gap is $1.417\times10^{-12}$. The discrepancy between the fine semidiscrete reference and the continuum solution is $2.842\times10^{-6}$ in the same relative RMS norm, so the target measures temporal error on this mesh. The selected runs include 8 stages in 0.427 s for TSFO-a, 8 stages in 0.428 s for TSFO-q, 12 stages in 0.215 s for Radau IIA(5). Radau IIA(5) has the lowest median time among the selected configurations. Each dash indicates that no sampled step count attained the target.
% END_AUTO_weno_budget_text
\FloatBarrier

\subsection{Measured work and robustness of the stage solves}\label{sec:cost_protocol}
\subsubsection{Timing at common endpoint error targets}
Times are medians of three single-thread integrations. The timed ODE runs return the endpoint. Separate untimed replays provide node histories for the transient-error diagnostics. Timings include residual and trajectory-rate evaluation, Jacobian assembly or actions, Newton updates, factorizations, and linear solves where required. The acoustic timings include Fourier preconditioner construction. The adaptive runs include their error control and factorization reuse. References, error evaluation, and plotting are outside the timed region.

At fixed step counts, \cref{tab:cost} separates the largest factored system from the number of factorizations. The TSFO Newton matrices have dimension $N$, and the directly coupled Gauss and three-stage Radau matrices have dimensions $2N$ and $3N$. SDIRK4 also uses matrices of dimension $N$, with more stages. The ROS4 implementation performs one factorization and four linear solves per step. Each TSFO Newton matrix also contains a derivative of the trajectory rate. These differences account for the distinct work associated with each stage count.
\begin{table}[htbp]\centering\small
\caption{Fixed-step ODE runs at 32 steps for HIRES and 100 for Robertson. Errors are endpoint infinity norms. LU is the number of factorizations, and the final column gives the largest factored dimension.}\label{tab:cost}
% BEGIN_AUTO_cost
\begin{tabular}{llrrrr}\toprule
Problem & Method & Endpoint error & Time (s) & LU & Max. system\\\midrule
HIRES & TSFO-a & $8.075\times10^{-10}$ & 0.0241 & 229 & 8\\
HIRES & TSFO-q & $2.022\times10^{-9}$ & 0.0264 & 230 & 8\\
HIRES & GL4 & $7.708\times10^{-9}$ & 0.0190 & 120 & 16\\
HIRES & Radau IIA(5) & $6.423\times10^{-11}$ & 0.0196 & 123 & 24\\
HIRES & SDIRK4 & $4.699\times10^{-9}$ & 0.0315 & 530 & 8\\
HIRES & ROS4 & $1.489\times10^{-7}$ & 0.0040 & 32 & 8\\
\midrule
Robertson & TSFO-a & $2.574\times10^{-12}$ & 0.0892 & 605 & 3\\
Robertson & TSFO-q & $6.395\times10^{-12}$ & 0.0529 & 605 & 3\\
Robertson & GL4 & $1.614\times10^{-12}$ & 0.0362 & 304 & 6\\
Robertson & Radau IIA(5) & $1.792\times10^{-12}$ & 0.0454 & 304 & 9\\
Robertson & SDIRK4 & $2.285\times10^{-11}$ & 0.0661 & 1222 & 3\\
Robertson & ROS4 & $1.755\times10^{-10}$ & 0.0097 & 100 & 3\\
\bottomrule\end{tabular}
% END_AUTO_cost
\end{table}

The work--precision comparison in \cref{tab:targets} uses the same endpoint infinity-norm error for all solvers on each problem and time interval. The table selects the fastest sampled configuration meeting the prescribed target for each method. The separately retained nodewise errors provide information on initial transients that may have decayed by the endpoint.
\begin{table}[htbp]\centering\small
\caption{Minimum sampled time in seconds for each prescribed endpoint error. The final times are $1$ for HIRES and van der Pol and $0.1$ for Robertson. A dash denotes an unmet target on the sampled configurations.}\label{tab:targets}
% BEGIN_AUTO_targets
\begin{tabular}{lrrr}\toprule
Method & HIRES, $10^{-7}$ & Robertson, $10^{-8}$ & van der Pol, $10^{-9}$\\\midrule
TSFO-a & 0.0128 & 0.0287 & 0.0143\\
TSFO-q & 0.0113 & 0.0275 & 0.0292\\
GL4 & 0.0190 & 0.0256 & 0.0129\\
Radau IIA(3) & 0.0494 & 0.0183 & 0.0102\\
Radau IIA(5) & 0.0051 & 0.0226 & 0.0129\\
SDIRK4 & 0.0162 & 0.0409 & 0.0470\\
ROS4 & 0.0107 & 0.0097 & 0.0061\\
HO4 & 0.0134 & 0.0153 & 0.0304\\
BDF adaptive & 0.0096 & 0.0082 & 0.0271\\
Radau adaptive & 0.0064 & 0.0036 & 0.0042\\
\bottomrule\end{tabular}
% END_AUTO_targets
\end{table}
% BEGIN_AUTO_targets_text
The lowest measured times among the sampled configurations are obtained by Radau IIA(5) for HIRES (0.0051 s), Radau adaptive for Robertson (0.0036 s), Radau adaptive for van der Pol (0.0042 s).
% END_AUTO_targets_text
\FloatBarrier

\subsubsection{Newton--Krylov robustness at larger amplitude}\label{sec:acoustic_pressure}
On the same acoustic grid we multiply both modal amplitudes by ten, to $0.02$ and $0.015$. The relaxation coefficient remains $4000$ and the final time remains $0.5$. We use $1$, $8$, $16$, and $32$ equal steps. The Newton stopping tolerance for the infinity norm of the preconditioned residual is $2\times10^{-11}\|u^0\|_\infty$. The Fourier preconditioner uses the linear-weight spatial symbol. Each GMRES solve uses a restart length of $24$ and a cap of five restart cycles. We attempt each configuration three times.
\begin{table}[H]\centering\small
\caption{Completion and work in the larger-amplitude acoustic test. Entries begin with completed runs out of three. For completed configurations we give the total GMRES iterations of the median-timed run and its time in seconds in parentheses.}\label{tab:acoustic_pressure}
% BEGIN_AUTO_acoustic_pressure
\begin{tabular}{rrccc}\toprule
Steps & $h$ & TSFO-a & TSFO-q & Radau IIA(5)\\\midrule
1 & 0.5 & 0/3, GMRES limit & 0/3, GMRES limit & 3/3, 24 (0.459)\\
8 & 0.0625 & 0/3, GMRES limit & 0/3, GMRES limit & 3/3, 32 (0.682)\\
16 & 0.03125 & 0/3, GMRES limit & 0/3, GMRES limit & 3/3, 40 (1.113)\\
32 & 0.015625 & 3/3, 189 (4.537) & 3/3, 197 (4.443) & 3/3, 70 (4.858)\\
\bottomrule\end{tabular}
% END_AUTO_acoustic_pressure
\end{table}
% BEGIN_AUTO_acoustic_pressure_text
The fewest sampled steps for which all repeats complete are 32 for TSFO-a, 32 for TSFO-q, and 1 for Radau IIA(5).
% END_AUTO_acoustic_pressure_text
The pressure test measures stage-solver completion at the sampled step sizes. The common-error comparison in \cref{tab:matched_pde} uses the smaller amplitudes.
\FloatBarrier

\section{Discussion, conclusions and future work}\label{sec:conclusion}
This work extends the explicit Li–Du method to a fully implicit, sequential two-stage, two-derivative time discretization. Direct Taylor matching yields fourth-order convergence at both the midpoint and the end of each step for fixed smooth systems. An analysis of the poles and the imaginary-axis boundary, together with the maximum modulus principle, characterizes the family of L-stable full-step updates. Within this family, \(C_a\) minimizes the magnitude of the leading classical error coefficient. The choice \(C_q=5/183\) of Huo and Su uniquely cancels the leading inverse-power term and gives quadratic stiff decay. Both choices use the same two sequential implicit solves.

The numerical experiments confirm fourth-order convergence on stiff ODEs, including Robertson kinetics, and on fixed spatial grids when the derivative of the discrete residual is evaluated consistently. In the fixed-step scalar test at \(z=-10^4\), the full-step amplification magnitude is \(6.26\times10^{-7}\) for TSFO-q and approximately \(0.999\) for the fourth-order Gauss–Legendre and symmetric Hermite–Obreshkov methods. The transport–relaxation experiments show that this damping advantage carries over to the fast component of the computed solution. At a common HIRES endpoint-error threshold of \(10^{-6}\), both TSFO variants require 16 implicit stages over the integration, compared with 24 for three-stage Radau IIA and 80 for SDIRK4. The sequential HIRES solves also reproduce the stacked Newton stage values to roundoff while reducing the dimension of each system from \(2N\) to \(N\). In the one-dimensional transport–relaxation test, TSFO-q reaches the prescribed endpoint accuracy faster than the sampled explicit Li–Du and adaptive Radau configurations.

The deep-stiff damping applies to the completed step. The midpoint can retain a fast component in the stiff limit, and the cost of evaluating the trajectory derivative and its Jacobian can offset the savings from smaller systems. For example, TSFO uses fewer stages in the two-dimensional acoustic comparison, yet three-stage Radau IIA has the shorter measured runtime. At larger perturbation amplitudes, the present Newton–Krylov implementation reaches its GMRES iteration limit at the coarser tested steps. The fixed-grid Burgers test also shows that an approximate physical interface rate can reduce the observed temporal order. These results motivate the development of less expensive interface rates consistent with the discrete residual, preconditioners tailored to the two stage matrices, and step-size control informed by stage convergence. A further theoretical question is under what conditions the fourth-order endpoint error bound can be made uniform in stiffness.

\paragraph{Data and code availability.}
The accompanying reproducibility package provides the implemented integrators, analytic derivatives, independent coefficient checks, the binary128 stiff ODE solver and its reference-resolution checks, all experiment scripts, recorded CSV/NPZ results, environment metadata, figure generation, and instructions for rebuilding the single-file LaTeX manuscript. The one- and two-dimensional spatial providers, stage-residual records and both-node refinement audit are included.

\paragraph{Acknowledgments.}
%This work is dedicated to the memory of my advisor, Professor Jiequan Li, whose encouragement motivated the explicit-to-implicit extension.
The author's research was supported by the Key Program of Department of Education of Henan Province (26A110007), the Young Talents Fund of Department of Science and Technology of Henan Province (252300423500), the Double First-Class Project of the School of Geomatics of Henan Polytechnic University (BSJJ202306), and the Doctoral Startup Foundation of Henan Polytechnic University (B2024-60).

\appendix
\section{Technical details for the global fourth-order proof}\label{app:convergence_details}
This appendix supplies the uniform bounds and full chain-rule calculations used in Section~\ref{sec:order}. Norms on matrices are induced by the fixed vector norm. Constants may change between estimates but are independent of $n,h$ on the fixed interval; they may depend on the system, $C,D$, and the chosen neighborhoods.

\subsection{Uniform local branches and stage sensitivities}\label{app:uniform_branches}
Let $\mathcal U_r=\{x:\operatorname{dist}(x,u([0,T]))\le r\}$, using the fixed vector norm. Choose $\rho>0$ so that the compact tube $\mathcal U_{3\rho}$ lies in the smoothness domain. For inputs $u^n\in\mathcal U_\rho$, use the closed stage ball $\overline B(u^n,\rho)\subset\mathcal U_{2\rho}$. All required derivative bounds are uniform on $\mathcal U_{3\rho}$. For a trial state $x$, the general stage residuals are
\begin{equation}\label{eq:app_general_residuals}
\begin{aligned}
 \mathcal R_{1/2}(x;u^n)
 &=x-u^n-\frac h4\bigl[\mathcal L(u^n)+\mathcal L(x)\bigr]
          -\frac{h^2}{48}\bigl[\mathcal H(u^n)-\mathcal H(x)\bigr],\\
 \mathcal R_1(x;u^n,u^{n+1/2})
 &=x-u^n-h\bigl[a_0\mathcal L(u^n)+a_{1/2}\mathcal L(u^{n+1/2})+a_1\mathcal L(x)\bigr]\\
 &\quad-h^2\bigl[b_0\mathcal H(u^n)+b_{1/2}\mathcal H(u^{n+1/2})+b_1\mathcal H(x)\bigr].
\end{aligned}
\end{equation}
Their derivatives with respect to $x$ are
\begin{equation}\label{eq:app_stage_jacobians}
\begin{aligned}
 M_{1/2}(x)&=I-\frac h4\mathcal L_u(x)+\frac{h^2}{48}\mathcal H_u(x),\\
 M_1(x)&=I-ha_1\mathcal L_u(x)-h^2b_1\mathcal H_u(x).
\end{aligned}
\end{equation}
In the following bounds, $\|\cdot\|_\infty$ denotes the supremum of the indicated operator norm on $\mathcal U_{2\rho}$. For fixed input $u^n$, consider $x\mapsto x-\mathcal R_{1/2}(x;u^n)$ and, with a midpoint value $y$ held fixed, $x\mapsto x-\mathcal R_1(x;u^n,y)$. Direct differentiation of \eqref{eq:app_general_residuals} gives
\begin{equation}\label{eq:app_fixedpoint_derivative_bound}
\begin{aligned}
 \|\partial_x[x-\mathcal R_{1/2}(x;u^n)]\|
 &\le\frac h4\|\mathcal L_u\|_\infty
       +\frac{h^2}{48}\|\mathcal H_u\|_\infty,\\
 \|\partial_x[x-\mathcal R_1(x;u^n,y)]\|
 &\le h|a_1|\|\mathcal L_u\|_\infty
       +h^2|b_1|\|\mathcal H_u\|_\infty.
\end{aligned}
\end{equation}
Both bounds are $O(h)$, uniformly in $y\in\overline B(u^n,\rho/4)$, and use no input derivative of $\Psi_h$. Bounded $\mathcal L,\mathcal H$ also make both right-hand sides $u^n+O(h)$. Choose a common $0<h_0\le1$ so they map $\overline B(u^n,\rho)$ into $\overline B(u^n,\rho/4)$ and contract by at most $1/2$. First apply the contraction mapping theorem to obtain $u^{n+1/2}=\Psi_h(u^n)$; then fix this value as $y$ and apply it to the endpoint. This gives unique roots in the selected stage balls, both $u^n+O(h)$. Shrink $h_0$ so exact trajectory increments over a step are also at most $\rho/4$.

For fixed known stage data, write $\mathcal R_k(x)=x+\mathcal E_k(x)$, where $\mathcal E_k$ has Lipschitz constant at most $K_{\mathrm{res}}h$. The triangle and reverse triangle inequalities give, on each selected ball,
\begin{equation}\label{eq:app_inverse_bound}
\begin{aligned}
 \|\mathcal R_k(x)-\mathcal R_k(y)\|
 &\le(1+K_{\mathrm{res}}h)\|x-y\|,\\
 \|\mathcal R_k(x)-\mathcal R_k(y)\|
 &\ge(1-K_{\mathrm{res}}h)\|x-y\|,
 \qquad k\in\{1/2,1\}.
\end{aligned}
\end{equation}
The contraction concerns $x\mapsto x-\mathcal R_k(x)$; \eqref{eq:app_inverse_bound} supplies continuity and inverse control of the residual. Here $K_{\mathrm{res}}$ is fixed independently of $h$. Shrink the common threshold so that $K_{\mathrm{res}}h_0\le1/2$. For every $0<h\le h_0$, the inverse residual bound is $(1-K_{\mathrm{res}}h)^{-1}=1+O(h)$. The same derivative estimates give $\|M_k(x)-I\|\le K_{\mathrm{res}}h$, so the Neumann series yields
\[
 \|M_k(x)^{-1}\|\le(1-K_{\mathrm{res}}h)^{-1}
                  \le1+2K_{\mathrm{res}}h\le2,
 \qquad M_k(x)^{-1}=I+O(h).
\]
The implicit function theorem now gives differentiability of the input-dependent local branches. Their uniqueness is restricted to the selected stage balls.

Here $\mathrm D$ denotes differentiation with respect to the input $u^n$, with $h$ fixed. Differentiating the first stage exactly gives
\begin{equation}\label{eq:app_mid_sensitivity}
 M_{1/2}(u^{n+1/2})\,\mathrm D\Psi_h(u^n)
 =I+\frac h4\mathcal L_u(u^n)+\frac{h^2}{48}\mathcal H_u(u^n).
\end{equation}
Differentiating the endpoint, including its dependence on $u^{n+1/2}=\Psi_h(u^n)$, gives
\begin{equation}\label{eq:app_end_sensitivity}
\begin{aligned}
 M_1(u^{n+1})\,\mathrm D\Phi_h(u^n)
 ={}&I+ha_0\mathcal L_u(u^n)+h^2b_0\mathcal H_u(u^n)\\
 &+\bigl[ha_{1/2}\mathcal L_u(u^{n+1/2})
          +h^2b_{1/2}\mathcal H_u(u^{n+1/2})\bigr]\mathrm D\Psi_h(u^n).
\end{aligned}
\end{equation}
The first identity implies $\mathrm D\Psi_h=I+O(h)$. Substitution in the second gives $\mathrm D\Phi_h=I+O(h)$: the initial-state terms have factors $h,h^2$, and the entire midpoint chain contribution is $O(h)$ times $I+O(h)$. These uniform identities establish both derivative norm bounds used in the main proof.

\emph{Closing the stepwise induction.} By the initial-error assumption, $\|e_0\|\le K_{\mathrm{init}}h^4$ for sufficiently small $h$. Shrink the common step threshold so that
\[
 e^{K_{\mathrm{prop}}T}
 \bigl(K_{\mathrm{init}}+K_{\mathrm{loc}}T\bigr)h_0^4<\rho/4.
\]
This choice also gives $\|e_0\|<\rho/4$, starting the induction. Assume $\|e_j\|<\rho/4$ through step $n$. Every segment joining $u^j$ to $u(t_j)$ then stays in $\mathcal U_\rho$, so the uniform estimates apply through step $n+1$. The recurrence gives $\|e_{n+1}\|<\rho/4$. Hence all steps with $t_n+h\le T$ exist, inputs stay in $\mathcal U_\rho$, and stages lie in $\mathcal U_{2\rho}$. For inexact solves with residuals at most $K_s h^5$, enlarge $K_{\mathrm{loc}}$ using \eqref{eq:app_residual_perturbations} before choosing $h_0$ and repeat the same induction. The common $h_0$ may depend on the fixed stiffness and spatial resolution through $\rho$ and the derivative bounds; no uniformity in stiffness or spatial resolution is asserted.

\subsection{Uniform local errors and the principal vector coefficient}\label{app:local_errors}
Since $u_t=\mathcal L(u)$ and $\mathcal L\in C^5$, repeated differentiation gives $u\in C^6$. Its time derivatives through degree six are bounded on $[0,T]$. Also $\mathcal H=\mathcal L_u\mathcal L\in C^4$. Taylor's formula with integral remainder therefore gives uniform remainders in the expansions used below.

The coefficient matches of Section~\ref{sec:method} imply
\begin{equation}\label{eq:app_exact_residual_bound}
\begin{aligned}
 &\|\mathcal R_{1/2}(u(t_n+h/2);u(t_n))\|\\
 &\quad+\|\mathcal R_1(u(t_n+h);u(t_n),u(t_n+h/2))\|\le Kh^5.
\end{aligned}
\end{equation}
Applying \eqref{eq:app_inverse_bound} first to the midpoint gives \eqref{eq:midlocal}. The difference between the endpoint residuals using $\Psi_h(u(t_n))$ and $u(t_n+h/2)$ has norm at most
\begin{equation}\label{eq:app_midpoint_substitution_bound}
\begin{aligned}
 &\bigl(h|a_{1/2}|\|\mathcal L_u\|_\infty
       +h^2|b_{1/2}|\|\mathcal H_u\|_\infty\bigr)
 \|\Psi_h(u(t_n))-u(t_n+h/2)\|\\
 &\qquad\le\bigl(h|a_{1/2}|\|\mathcal L_u\|_\infty
       +h^2|b_{1/2}|\|\mathcal H_u\|_\infty\bigr)
       K_{\mathrm{loc}}h^5=O(h^6).
\end{aligned}
\end{equation}
The $\mathcal L$ and $\mathcal H$ contributions are $O(h^6)$ and $O(h^7)$, respectively. The endpoint inverse estimate proves the uniform upper bound \eqref{eq:local}.

For the sharper asymptotic coefficient, use $\mathcal L(u(t))=u_t(t)$ and $\mathcal H(u(t))=u_{tt}(t)$ along the vector trajectory. Taylor expansion of the endpoint quadrature gives
\begin{equation}\label{eq:app_vector_quadrature}
\begin{aligned}
 &h\sum_{\theta\in\{0,1/2,1\}}a_\theta u_t(t_n+\theta h)
 +h^2\sum_{\theta\in\{0,1/2,1\}}b_\theta u_{tt}(t_n+\theta h)
 -\bigl[u(t_n+h)-u(t_n)\bigr]\\
 &\quad=h^5\left(\frac{a_{1/2}}{384}+\frac{a_1}{24}
                 +\frac{b_{1/2}}{48}+\frac{b_1}{6}-\frac1{120}\right)u^{(5)}(t_n)+O(h^6).
\end{aligned}
\end{equation}
All lower powers vanish by \eqref{eq:secondconditions}. Substitution of \eqref{eq:coefficients} evaluates the parenthesis as
\begin{equation}\label{eq:app_principal_coefficient}
 \frac{1-60C+15D}{2880},
 \qquad D=-C\ \Longrightarrow\ \frac{1-75C}{2880}.
\end{equation}
The left side of \eqref{eq:app_vector_quadrature} is the negative exact-endpoint residual. Replacing its midpoint evaluations $\mathcal L(u(t_n+h/2))$ and $\mathcal H(u(t_n+h/2))$ by $\mathcal L(\Psi_h(u(t_n)))$ and $\mathcal H(\Psi_h(u(t_n)))$ changes this left-hand side by $O(h^6)$, as \eqref{eq:app_midpoint_substitution_bound} shows. Its $h^5$ coefficient is therefore unchanged. The segment-averaged endpoint Jacobian is $I+O(h)$, so inversion preserves that coefficient and changes only the $O(h^6)$ remainder. This proves \eqref{eq:principal} for the vector problem directly. The uniform upper bounds \eqref{eq:midlocal}--\eqref{eq:local} and this sharper expansion are separate conclusions of the same Taylor calculation.

\subsection{Inexact solves and local Newton convergence}\label{app:inexact_newton}
For fixed input $u^n$, let tildes denote approximate stage solves in the same local neighborhoods. The inverse bound and the midpoint dependence of the second residual give
\begin{equation}\label{eq:app_residual_perturbations}
\begin{aligned}
 \|\widetilde u^{n+1/2}-\Psi_h(u^n)\|
 &\le K\|\mathcal R_{1/2}(\widetilde u^{n+1/2};u^n)\|,\\
 \|\widetilde u^{n+1}-\Phi_h(u^n)\|
 &\le K\bigl(\|\mathcal R_1(\widetilde u^{n+1};u^n,\widetilde u^{n+1/2})\|\\
 &\hspace{20mm}+h\|\mathcal R_{1/2}(\widetilde u^{n+1/2};u^n)\|\bigr).
\end{aligned}
\end{equation}
Thus residuals at most $K_s h^5$ add at most $Kh^5$ to the endpoint recurrence, as required in \cref{thm:convergence}. More generally, if every unscaled stage residual is bounded by $\tau_h$ and the computed states stay in the selected neighborhoods, the same summation gives
\begin{equation}\label{eq:app_inexact_global_bound}
\begin{aligned}
 &\max_{0\le nh\le T}\|u^n-u(t_n)\|
 +\max_{t_n+h\le T}\|u^{n+1/2}-u(t_n+h/2)\|\\
 &\qquad\le K_T\left(\|u^0-u(0)\|+h^4+\frac{\tau_h}{h}\right).
\end{aligned}
\end{equation}
This makes the sufficient choice $\tau_h=O(h^5)$ explicit. A fixed floating-point correction tolerance alone does not enforce this choice as $h\to0$; roundoff and the attained final residual must be distinguished from the exact-arithmetic convergence statement.

For Newton's method, hold the known data of either stage fixed. Since $\mathcal L\in C^5$ and $\mathcal H\in C^4$, both $\mathcal L_{uu}$ and $\mathcal H_{uu}$ are bounded; the latter includes third state derivatives of $\mathcal L$. In particular,
\begin{equation}\label{eq:app_newton_matrix_derivative_bound}
\begin{aligned}
 \|\mathrm D M_{1/2}(x)\|
 &\le\frac h4\|\mathcal L_{uu}\|_\infty
       +\frac{h^2}{48}\|\mathcal H_{uu}\|_\infty
       \le C_{\mathrm{der}}h,\\
 \|\mathrm D M_1(x)\|
 &\le h|a_1|\|\mathcal L_{uu}\|_\infty
       +h^2|b_1|\|\mathcal H_{uu}\|_\infty
       \le C_{\mathrm{der}}h.
\end{aligned}
\end{equation}
Equation~\eqref{eq:app_newton_matrix_derivative_bound} gives $\|M_k(x)-M_k(y)\|\le C_{\mathrm{der}}h\|x-y\|$ on a stage ball, while \eqref{eq:app_inverse_bound} gives a separate uniform inverse bound $C_{\mathrm{inv}}\le2$. For a local root $x_*$, a full Newton update obeys the exact integral identity
\[
\begin{aligned}
 x^{(\ell+1)}-x_*
 =M_k(x^{(\ell)})^{-1}\int_0^1
 &\bigl[M_k(x^{(\ell)})
       -M_k(x_*+s(x^{(\ell)}-x_*))\bigr]\\
 &\qquad\cdot(x^{(\ell)}-x_*)\,ds.
\end{aligned}
\]
Consequently,
\begin{equation}\label{eq:app_newton_local_bound}
 \|x^{(\ell+1)}-x_*\|
 \le C_{\mathrm N}h\|x^{(\ell)}-x_*\|^2,
 \qquad C_{\mathrm N}=\tfrac12C_{\mathrm{inv}}C_{\mathrm{der}}.
\end{equation}
All constants are independent of $n,h$ for $0<h\le h_0$. Set $r_{\mathrm N}=\rho/2$. Since $\|x_*-u^n\|\le\rho/4$, the triangle inequality gives $\overline B(x_*,r_{\mathrm N})\subset\overline B(u^n,3\rho/4)\subset\overline B(u^n,\rho)$. Thus the Newton ball lies inside the selected stage ball. Shrink $h_0$ so $C_{\mathrm N}h_0r_{\mathrm N}\le1/2$. Then \eqref{eq:app_newton_local_bound} keeps any Newton iterate initially in this smaller ball inside it and gives quadratic convergence. The predictors $u^n$ for the first solve and $u^{n+1/2}$ for the second are $O(h)$ from their respective roots, so enter this region for sufficiently small $h$. This is a local result; Sections~\ref{sec:newton} and \ref{sec:cost_protocol} describe backtracking and failure reporting for larger steps. The recorded fixed-tolerance experiments concern their resolved refinement range.

\subsection{A polynomial reaction example}\label{app:robertson_regular}
The vector field in \eqref{eq:rober} conserves $u_1+u_2+u_3$. On the respective faces $u_i=0$ of
\[
 \Delta=\{u\in\R^3:u_i\ge0,\ \textstyle\sum_i u_i=1\},
\]
the zero-component derivatives are $10^4u_2u_3$, $0.04u_1$ and $3\times10^7u_2^2$, all nonnegative. To handle tangency at the boundary, use $\mathcal L+\varepsilon(\mathbf1-3u)$, $0<\varepsilon\le1$, with the same initial value in $\Delta$. The sum of its components is $3\varepsilon(1-\sum_i u_i)$, so it preserves the plane $\sum_i u_i=1$. On a face $u_i=0$,
\[
 [\mathcal L+\varepsilon(\mathbf1-3u)]_i
 =\mathcal L_i+\varepsilon\ge\varepsilon>0.
\]
It is therefore strictly inward on every face. A first-exit argument keeps each perturbed solution in $\Delta$, giving existence on every finite interval. The perturbed fields have uniform bounds and Lipschitz constants near $\Delta$; Gronwall's inequality gives
\[
 \sup_{0\le t\le T}\|u_\varepsilon(t)-u_\delta(t)\|
 \le C_T|\varepsilon-\delta|,\qquad 0<\varepsilon,\delta\le1.
\]
Thus the solutions converge uniformly on every fixed finite interval as $\varepsilon\downarrow0$. Passing to the limit in the integral equation gives the original solution, which remains in $\Delta$ and exists for every finite $T$.

The field is polynomial, so it and all state derivatives needed in Section~\ref{sec:order} are bounded on a neighborhood of $\Delta$. Thus \cref{thm:convergence} applies on every fixed finite interval for any fixed $C,D$ in \eqref{eq:coefficients}. The constants can depend on the fixed reaction rates, parameters and $T$.

\end{document}